\documentclass[a4paper,11pt]{article}
\usepackage{amsfonts,amsmath,amsthm}
\usepackage{graphicx}
\usepackage{color}

\newcommand{\A}{\mathcal{A}}
\newcommand{\E}{\mathbb{E}}
\newcommand{\erre}{\mathbb{R}}

\newcommand{\cp}[2]{\langle#1,#2\rangle}
\newcommand{\ds}{\displaystyle}
\newcommand{\tr}{\mathop{\mathrm{Tr}}\nolimits}

\newtheorem{prop}{Proposition}[section]
\newtheorem{thm}[prop]{Theorem}
\newtheorem{coroll}[prop]{Corollary}
\newtheorem{lemma}[prop]{Lemma}
\theoremstyle{definition}
\newtheorem{rmk}[prop]{Remark}

\newsymbol\lesssim 132E

\title{On controlled linear diffusions with delay in a model of
  optimal advertising under uncertainty with memory effects}

\author{Fausto Gozzi\thanks{Facolt\`a di Economia, Libera Universit\`a
    degli Studi Sociali ``Guido Carli'', 00198 Roma, Italy. e-mail:
    \texttt{fgozzi@luiss.it}.} \and Carlo Marinelli\thanks{Institut
    f\"ur Angewandte Mathematik, Universit\"at Bonn, Wegelerstr. 6,
    D-53115 Bonn, Germany. URL:
    \texttt{http://www.uni-bonn.de/$\sim$cm788}.} \and Sergei
  Savin\thanks{Graduate School of Business, Columbia University, New
    York, NY 10027, USA. e-mail: \texttt{svs30@columbia.edu}.}}

\date{\normalsize June 7, 2007. Revised July 21, 2008.}

\begin{document}
\maketitle

\begin{abstract}
  We consider a class of dynamic advertising problems under
  uncertainty in the presence of carryover and distributed forgetting
  effects, generalizing the classical model of Nerlove and Arrow
  \cite{NA}. In particular, we allow the dynamics of the product
  goodwill to depend on its past values, as well as previous
  advertising levels. Building on previous work (\cite{levico}), the
  optimal advertising model is formulated as an infinite dimensional
  stochastic control problem. We obtain (partial) regularity as well
  as approximation results for the corresponding value function. Under
  specific structural assumptions we study the effects of delays on
  the value function and optimal strategy. In the absence of carryover
  effects, since the value function and the optimal advertising policy
  can be characterized in terms of the solution of the associated HJB
  equation, we obtain sharper characterizations of the optimal policy.
  
\end{abstract}

\section{Introduction}
This paper is devoted to the study of a class of optimal control
problems for linear stochastic differential equations with delay both
in the state and the control term, and is a natural continuation of
\cite{levico}. These problems arise in the theory of optimal
advertising under uncertainty with memory structures. We approach the
problem using stochastic control techniques in infinite dimensions.

In particular, in \cite{levico} we considered a controlled stochastic
differential equation (SDE) with delay entering both the state and the
control variable as an extension of the dynamic advertising model of
Nerlove and Arrow \cite{NA}. The results of \cite{levico} are the
following: we construct a controlled infinite dimensional SDE that is
equivalent to the controlled SDE with delay, we prove a verification
theorem, and we exhibit a simple example for which the Bellman
equation associated to the control problem admits a sufficiently
regular solution, hence the verification theorem can be applied. 
In the present manuscript, we extend \cite{levico} by
  developing several new sets of results.  On the one hand, we provide
  qualitative characterization of the first- and second-order
  properties of the optimal value function (Section
  \ref{sec:reformul}.1). In particular, we show that, under natural
  restrictions, the monotonicity of the optimal value function with
  respect to the initial goodwill profile still holds even in the
  presence of the state and control-related delay terms in the
  advertising dynamics (Proposition \ref{prop:Vincr}). In addition, we
  establish that the decreasing marginal influence of the attained
  goodwill levels on the primal profit components is retained by the
  optimal profit function (Proposition \ref{prop:Vconv}). As is well
  known, this last property is important in reducing the computational
  load required to solve the time- and space-discretized version of
  our problem by dynamic programming methods.
  On the other hand, in the view of intractability of the general
  variant of our problem, we propose approximation schemes for the
  optimal value function (Theorem \ref{prop:approx}) and for the
  optimal advertising policy (Propositions 3.11 and 3.12). The latter
  result is of particular importance since it suggests a
  computationally feasible approach to constructing asymptotically
  optimal advertising trajectories.
  In addition, we provide a complete characterization of the optimal
  advertising policy in the case when the cost function is quadratic
  and the reward function is linear in goodwill level (Section
  3.3). For a specific instance of this case we conduct a numerical
  study aimed at demonstrating the importance of proper accounting of
  the delay effects in calculating the optimal advertising policy.

  Finally, we are able to provide sharper characterization of the
  optimal policies in the case when the influence of advertising on
  the goodwill evolution is instantaneous and the delay effects are of
  the state-only type (Section \ref{sec:ex}). The key result in this
  section is Theorem \ref{thm:L2} which formulates sufficient
  conditions ensuring that the optimal advertising policy is of the
  feedback type. In particular, in the case of linear cost function
  the optimal control takes a particularly simple ``bang-bang'' form
  (Corollary \ref{cor:bang}).


Optimal control problems for stochastic systems with delay in the
state term admit alternative, more traditional treatments: for
instance, see \cite{elsanosi} and \cite{larssen} for a more direct
application of the dynamic programming principle without appealing to
infinite dimensional analysis, and \cite{KSlibro} for the
linear-quadratic case. However, we would like to point out that none
of the methods just mentioned apply to the control of stochastic
differential equations with delay in the control term.

Analysis of advertising policies has always been occupying a
front-and-center place in the marketing research. The sheer size of
the advertising market (over \$143 billion in the US in 2005
\cite{admarket}) and the strong body of evidence of systematic
over-advertising by firms across many industries (see e.g.
\cite{overad1}, \cite{overad4}, \cite{overad2}, \cite{overad3})
has caused a renewal of attention to the proper accounting for the
so-called ``carryover'' or ``distributed lag'' advertising effects.
The term ``carryover'' designates an empirically observed
advertising feature under which the advertising influence on product
sales or goodwill level is not immediate, but rather is spread over
some period of time: according to a survey of recent empirical
``carryover'' research by Leone \cite{leone}, delayed advertising effects
can last between 6 and 9 months in different settings.

On the theoretical front, pioneering work of \cite{VW} and \cite{NA}
has paved the way for the development of a number of models dealing
with the optimal distribution of advertising spending over time in
both monopolistic and competitive settings. A comprehensive review of
the state of the advertising control literature in \cite{sethi} points
out that the majority of these models operate under deterministic
assumptions and do not capture some of the most essential
characteristics of real-world advertising phenomena. On the empirical
side, one of the first and most important substreams of advertising
literature was formed by the papers focused on the studies of
distributed advertising lag (see e.g. \cite{bass2},
\cite{bassparsons}, \cite{griliches}). An important early empirical
result was obtained by Bass and Clark \cite{bassclark}, who
established that the initially adopted models with monotone decreasing
lags (see \cite{Koyck}) are often inferior in their explanatory power
to the models with more general lag distributions.

Despite the wide and growing empirical literature on the measurement
of carryover effects, there are practically no analytical studies that
incorporate distributed lag structure into the optimal advertising
modeling framework in the stochastic setting. The only papers dealing
with optimal dynamic advertising with distributed lags we are aware of
are \cite{bassparsons} (which provides a numerical solution to a
discrete-time deterministic example), and \cite{hartl}, \cite{HS}
(which applies a version of the maximum principle in the deterministic
setting). The creation of models which incorporate the treatment of
``carryover'' effects in the stochastic settings have long been
advocated in the advertising modeling literature (see e.g.
\cite{hartl}, \cite{sethi} and references therein).

As mentioned above, in this work we study a class of stochastic models
deriving from that of Nerlove and Arrow \cite{NA}, incorporating both
the advertising lags as well as distributed ``churn''
(\cite{sethi-VWcomp}), or ``forgetting'', effects.  More precisely,
we formulate an optimization program that seeks to maximize the
goodwill level at a given time $T>0$ net of the cumulative cost of
advertising until $T$. This optimization problem is studied using
techniques of stochastic optimal control in infinite dimensions, using
the modeling approach of \cite{levico}: in particular, we specify the
goodwill dynamics in terms of a controlled stochastic delay
differential equation (SDDE), that can be rewritten as a stochastic
differential equation (without delay) in a suitable Hilbert space.
This allows us to associate to the original control problem for the
SDDE an equivalent infinite dimensional control problem for the
``lifted'' stochastic equation.

The paper is organized as follows: in section \ref{sec:formul} we
formulate the optimal advertising problem as an optimal control
problem for an SDE with delay, and we recall the equivalence result of
\cite{levico}. In section \ref{sec:reformul} we prove the above
mentioned results about the value function and approximate strategies
in the general case, together with a detailed discussion of the effect
of delays in a specific situation. Section \ref{sec:ex} treats the
case of distributed forgetting in the absence of advertising
carryover.

Let us conclude this introduction fixing notation and recalling some
notions that will be needed. Given a lower semicontinuous convex
function $f:E\to\bar{\erre}:=\erre\cup\{+\infty\}$ on a Hilbert space
$E$ with inner product $\cp{\cdot}{\cdot}$, we denote its conjugate by
$f^*(y):=\sup_{x\in E}(\cp{x}{y}-f(x))$. Recall also that $D^-
f^*(y)=\arg\max_{x\in E} (\cp{x}{y}-f(x))$, where $D^-$ stands for the
subdifferential operator (see e.g. \cite{barbu-v}, p.~103).
Throughout the paper, $X$ will be the Hilbert space defined as
$$
X = \erre \times L^2([-r,0],\erre),
$$
with inner product
$$
\langle x,y\rangle = x_0y_0 + \int_{-r}^0 x_1(\xi)y_1(\xi)\,d\xi
$$
and norm
$$
|x| = \left( |x_0|^2 + \int_{-r}^0 |x_1(\xi)|^2\,d\xi\right)^{1/2},
$$
where $r>0$, $x_0$ and $x_1(\cdot)$ denote the $\erre$-valued and the
$L^2([-r,0],\erre)$-valued components, respectively, of the generic
element $x$ of $X$. Given $f:X\to\erre$, $k\in\{0,1\}$, we shall
denote by $\partial_k f$, $D^-_k f$, respectively, the partial
derivative and the subdifferential of $f$ with respect to the $k$-th
component.  We shall use mollifiers in a standard way: for $\zeta \in
C^\infty(\erre^d,\erre_+)$, equal to zero for $|x|>1$ and such that
$\int_{\erre^d}\zeta(x)\,dx=1$, we shall set
$\zeta_\lambda(x)=\lambda^{-d}\zeta(\lambda^{-1}x)$ for $\lambda\neq 0$.
By $a \lesssim b$ we mean that there exists a constant $N$ such that
$a \leq Nb$. If $N$ depends on some parameter of interest $p$, we
shall write $N(p)$ and $a \lesssim_p b$.

\section{The model}\label{sec:formul}
We consider a monopolistic firm preparing the market introduction of a
new product at some time $T$ in the future. In defining the state
descriptor for a firm to follow we use the Nerlove-Arrow framework and
consider the product's ``goodwill stock'' $y(t)$, $0\leq s \leq t \leq
T$. The firm directly influences the rate of advertising spending
$z(t)$ to induce the following trajectory for the goodwill stock:
\begin{equation}
\label{eq:SDDE}
\left\{\begin{array}{l}
dy(t) = \ds \left[a_0 y(t) + \int_{-r}^0 a_1(\xi)y(t+\xi)\,d\xi
        + b_0 z(t) +
            \int_{-r}^0b_1(\xi)z(t+\xi)\,d\xi\right]dt \\[10pt]
\ds \qquad\qquad  + \sigma\, dW_0(t), \quad s \leq t\leq T \\[10pt]
y(s)=x_0; \quad y(s+\xi)=x_1(\xi),\;
z(s+\xi)=\delta(\xi),\;\; \xi\in[-r,0],
\end{array}\right.
\end{equation}
where the Brownian motion $W_0$ is defined on a filtered probability
space $(\Omega,\mathcal{F},\mathbb{F}=(\mathcal{F}_t)_{t\geq 0},%
\mathbb{P})$, with $\mathbb{F}$ being the completion of the filtration
generated by $W_0$. We assume that the advertising spending rate
$z(t)$ is constrained to remain in the set $\mathcal{U}$, the space of
$\mathbb{F}$-adapted processes taking values in a compact
  interval $U\subseteq\erre_+$. In addition, we assume that the
following conditions hold:
\begin{itemize}
\item[(i)] $a_0 \leq 0$;
\item[(ii)] $a_1(\cdot) \in L^2([-r,0],\erre)$;
\item[(iii)] $b_0 \geq 0$;
\item[(iv)] $b_1(\cdot) \in L^2([-r,0],\erre_+)$;\label{iv}
\item[(v)] $x_0 \geq 0$;
\item[(vi)] $x_1(\cdot) \geq 0$, with $x_1(0)=x_0$;
\item[(vii)] $\delta(\cdot) \geq 0$.
\end{itemize}
Here $a_0$ and $a_1(\cdot)$ describe the process of goodwill
deterioration when the advertising stops, and $b_0$ and $b_1(\cdot)$
provide the characterization of the effect of the current and the past
advertising rates on the goodwill level. The values of $x_0$,
$x_1(\cdot)$ and $\delta(\cdot)$ reflect the ``initial'' goodwill and
advertising trajectories.  Note that we recover the model of Nerlove
and Arrow from (\ref{eq:SDDE}) in the deterministic setting
($\sigma=0$) in the absence of delay effects
($a_1(\cdot)=b_1(\cdot)=0$).

Setting $X \ni x:=(x_0,x_1(\cdot))$ and denoting by $y^{s,x,z}(t)$,
$t\in[0,T]$, a solution of (\ref{eq:SDDE}), we define the objective
functional
\begin{equation}
\label{eq:obj-orig}
J(s,x;z) =
\E\left[\varphi_0(y^{s,x,z}(T))-\int_s^T h_0(z(t))\,dt \right],
\end{equation}
where $\varphi_0:\erre\to\erre$ and $h_0:\erre_+\to\erre_+$ are
measurable utility and cost functions, respectively, satisfying the
conditions
\begin{equation}\label{eq:hphipol}
|\varphi_0(x)| \leq K(1+|x|)^m
\qquad \forall x\in\erre,
\end{equation}
and
\begin{equation}
  \label{eq:schip}
  |h(x)| \leq K \qquad \forall x \in U,
\end{equation}
for some $K>0$ and $m\geq 0$. In the sequel we shall often move the
superscripts $s,x,z$ to the expectation sign, with obvious meaning of
the notation.
Let us also define the value function $V$ for this problem as follows:
$$
V(s,x) = \sup_{z\in\mathcal{U}} J(s,x;z).
$$
We shall say that $z^*\in\mathcal{U}$ is an optimal strategy if it
is such that
$$
V(s,x)=J(s,x;z^*).
$$
The problems we will deal with are the maximization of the objective
functional $J$ over all admissible strategies $\mathcal{U}$, and the
characterization of the value function $V$ and of the optimal strategy
$z^*$.

Throughout the paper we will always assume that the assumptions of
this section hold true. In particular the constants $T$, $m$ and $K$
are fixed from now on.

\subsection{An equivalent infinite dimensional Markovian
  representation}
We shall recall a representation result (proposition \ref{prop:equiv}
below) proved in \cite{levico}, generalizing a corresponding
deterministic result due to Vinter and Kwong \cite{VK}.

\noindent
Let us define an operator $A:D(A)\subset X\to X$ as
follows:
\begin{eqnarray*}
A: (x_0,x_1(\xi)) &\mapsto& \Big(a_0x_0 + x_1(0),a_1(\xi)x_0
                           -{dx_1(\xi)\over d\xi}\Big) \quad%
\textrm{a.e.}\ \xi\in [-r,0], \\
D(A) &=& \left\{ x \in X:
       x_1 \in W^{1,2}([-r,0];\erre),\; x_1(-r)=0\right\}.
\end{eqnarray*}
Moreover, define the bounded linear control operator $B:U \to X$ as
\begin{equation}
B: u \mapsto \Big(b_0u,b_1(\cdot)u\Big),
\end{equation}
and finally the operator $G:\erre \to X$ as $G: x_0 \mapsto (\sigma
x_0, 0)$. Sometimes it will be useful to identify the operator $B$
with the element $(b_0,b_1)\in X$.
\begin{prop}\label{prop:equiv}
  Let $Y(\cdot)$ be the weak solution of the abstract evolution
  equation
\begin{equation}
\label{eq:abstract}
\left\{\begin{array}{l}
dY(t) = (A Y(t)+Bz(t))\,dt + G\,dW_0(t) \\[8pt]
Y(s) = \bar{x} \in X,
\end{array}\right.
\end{equation}
with arbitrary initial datum $\bar{x} \in X$ and control
$z\in\mathcal{U}$. Then, for $t\geq r$, one has,
$\mathbb{P}$-a.s.,
$$
Y(t) = M(Y_0(t),Y_0(t+\cdot),z(t+\cdot)),
$$
where
\begin{eqnarray*}
M: X \times L^2([-r,0],\erre) &\to& X \\
(x_0,x_1(\cdot),v(\cdot)) &\mapsto& (x_0,m(\cdot)),
\end{eqnarray*}
$$
m(\xi) := \int_{-r}^\xi a_1(\zeta) x_1(\zeta-\xi)\,d\zeta
+ \int_{-r}^\xi b_1(\zeta) v(\zeta-\xi)\,d\zeta.
$$
Moreover, let $\{y(t),\; t\geq -r\}$ be a continuous solution of the
stochastic delay differential equation (\ref{eq:SDDE}), and $Y(\cdot)$
be the weak solution of the abstract evolution equation
(\ref{eq:abstract}) with initial condition
$$
\bar{x} = M(x_0,x_1,\delta(\cdot)).
$$
Then, for $t\geq 0$, one has, $\mathbb{P}$-a.s.,
$$
Y(t) = M(y(t),y(t+\cdot),z(t+\cdot)),
$$
hence $y(t)=Y_0(t)$, $\mathbb{P}$-a.s., for all $t\geq 0$.
\end{prop}
Using this equivalence result, we can now give a Markovian
reformulation on the Hilbert space $X$ of the problem of maximizing
(\ref{eq:obj-orig}), as in \cite{levico}.
In particular, denoting by $Y^{s,\bar{x},z}(\cdot)$ a mild solution of
(\ref{eq:abstract}), (\ref{eq:obj-orig}) is equivalent to
\begin{equation}
\label{ex1jbis}
J(s,x;z) =
\E\left[ \varphi(Y^{s,\bar{x},z}(T))+\int_s^T h(z(t))\,dt \right],
\end{equation}
with the functions $h:U\to\erre$ and $\varphi:X\to\erre$ defined by
\begin{eqnarray*}
h(z) &=& -h_0(z) \\
\varphi(x_0,x_1) &=& \varphi_0(x_0).
\end{eqnarray*}
Hence also $V(s,x)=\sup_{z\in \mathcal{U}} J(s,x;z)$.

\section{The case of delay in the state and the control term}
\label{sec:reformul}
The aims of this section are the following: to prove regularity
properties of the value function, to develop an approximation scheme
for the value function and the optimal strategy, and to illustrate in
a numerical example the effects of the delay structures in our model.
In particular, we prove that, under natural assumptions, the value
function is continuous in both arguments, and monotone concave with
respect to the initial goodwill profile. As already remarked, this
property is essential in order to obtain computationally tractable
discrete-time and discrete-state-space dynamic programming versions of
our problem.

Moreover, since we cannot guarantee that the Bellman equation
associated to our control problem admits a solution in general (nor do
we have any information about its uniqueness and regularity), it is of
primary interest to obtain approximation schemes for the optimal value
function and for the optimal advertising policy. The latter result is
of particular importance since it suggests a computationally feasible
approach to constructing asymptotically optimal advertising
trajectories.

In addition, in the last subsection we provide a complete
characterization of the optimal advertising policy in the case when
the cost function is quadratic and the reward function is linear in
goodwill level, and for a specific instance of this case we
conduct a numerical study aimed at demonstrating the importance of
proper accounting of the delay effects in calculating the optimal
advertising policy.

\subsection{Qualitative properties of the value function}
Let us first show that the value function is finite.
\begin{prop}\label{prop:Vfin}
  There exists a
  constant $N=N(T,m,K)$ such that $|V(s,x) |\leq N(1+|x|)^m$ for all
  $s\in[0,T]$, $x\in X$.
\end{prop}
\begin{proof}
  The estimate from below simply follows by taking a constant
  deterministic control. For the estimate from above we have,
  recalling that $h(x) \leq 0$ for all $x\in U$,
  \begin{eqnarray*}
    V(s,x) &\leq& \sup_{z\in\mathcal{U}} \E_{s,x}^z \Big[
    \int_s^T h(z(t))\,dt + |\varphi(Y(T))| \Big]\\
    &\leq& K \sup_{z\in\mathcal{U}} \E_{s,x}^z (1+|Y(T)|)^m.
  \end{eqnarray*}
  Moreover, we have
  \begin{align*}
    \E|Y(T)|^m &\lesssim |e^{(T-s)A}x|^m 
         + \E\bigg|\int_0^{T-s} e^{(T-s-r)A}Bz(r)\,dr\bigg|^m
    + \E\bigg|\int_0^{T-s} e^{(T-s-r)A}G\,dW(r)\bigg|^m\\
    &\lesssim_T M_T^m |x|^m + |B|^mM_T^m \bar{z}^m + \E|W_A(T)|^m,
  \end{align*}
  where $\bar{z}:=\max\{z:\,z\in U\}$, $W_A(t):=\int_0^{t-s}
  e^{(t-s-r)A}G\,dW(r)$, and
  $M_T:=\sup_{t\in[0,T]}|e^{tA}|$. Recalling that $GG^*$ is of trace
  class, hence
  \[
  \E|W_A(T)|^2 = \mathrm{Tr}\,\int_0^{T-s} e^{tA}GG^*e^{tA^*}\,dt <
  \infty
  \]
  (see e.g. \cite{DP-K}, Proposition 2.2), i.e. $W_A(T)$ is a
  well-defined Gaussian random variable on $X$, we also get that
  $\E|W_A(T)|^m<\infty$. The proof is completed observing that the
  upper bound on $\E|Y(T)|^m$ is uniform over $z$.
\end{proof}
%

We establish now some qualitative properties of the value function
that do not require studying an associated Bellman equation. The
following simple result, typical of control problems with linear
dynamics, asserts that the value function inherits the concavity with
respect to the space variable from the reward and cost functions.
\begin{prop}\label{prop:Vconv}
  If $\varphi$ and $h$ are concave,
  then the value function $V(s,x)$ is proper concave with respect to
  $x$.
\end{prop}
\begin{proof}
Properness follows by the previous proposition. Moreover, let $x^1$,
$x^2 \in X$. Then
\begin{eqnarray*}
\lambda V(s,x^1) + (1-\lambda)V(s,x^2) &=&
\lambda \sup_{z\in\mathcal{U}}
\E\Big[\int_s^T h(z(t))\,dt + \varphi(y^{s,x^1,z}(T))\Big] \\
&& + (1-\lambda) \sup_{z\in\mathcal{U}}
\E\left[\int_s^T h(z(t))\,dt + \varphi(y^{s,x^2,z}(T))\right] \\
&=& \sup_{z^1,z^2\in\mathcal{U}}
\E\bigg[\int_s^T [\lambda h(z^1(t))+(1-\lambda)h(z^2(t))]\,dt \\
&& \phantom{\sup_{z^1,z^2\in\mathcal{Z}}\E\bigg[}
+\lambda\varphi(y^{s,x^1,z^1}(T)) + (1-\lambda)\varphi(y^{s,x^2,z^2}(T))\bigg].
\end{eqnarray*}
Since $\mathcal{U}$ is a convex set and $h$ is concave, then
$z_\lambda:=\lambda z^1 + (1-\lambda) z^2$ is admissible for any
choice of $z^1$, $z^2\in\mathcal{U}$, and one has
\begin{equation}
  \label{eq:hc}
  h(z_\lambda(s)) \geq \lambda h(z^1(s))+(1-\lambda)h(z^2(s)).
\end{equation}
Moreover, by linearity
of the state equation, it is easy to prove that
$$
y^{s,x_\lambda,z_\lambda}(T) = \lambda y^{s,x^1,z^1}(T) +
(1-\lambda) y^{s,x^2,z^2}(T),
$$
hence, by the concavity of $\varphi$,
\begin{equation}
  \label{eq:phic}
  \varphi(y^{s,x_\lambda,z_\lambda}(T)) \geq
  \lambda \varphi(y^{s,x^1,z^1}(T)) +
  (1-\lambda) \varphi(y^{s,x^2,z^2}(T)).
\end{equation}
Therefore, as a consequence of (\ref{eq:hc}) and (\ref{eq:phic}), we obtain
\begin{eqnarray*}
\lambda V(s,x^1) + (1-\lambda)V(s,x^2)
&\leq&
\sup_{z_\lambda\in\mathcal{U}} \E\bigg[\int_s^T h(z_\lambda(t))\,dt +
\varphi(y^{s,x_\lambda,z_\lambda}(T)) \bigg] \\
&\leq&
\sup_{z\in\mathcal{U}} \E\bigg[\int_s^T h(z(t))\,dt +
\varphi(y^{s,x_\lambda,z}(T)) \bigg] \\
&=& V(s,\lambda x^1+(1-\lambda)x^2),
\end{eqnarray*}
which proves the claim.
\end{proof}
As a consequence of the previous propositions we obtain the following
regularity result. Of course it would be ideal to obtain a result
guaranteeing that $V \in C^{0,1}([0,T]\times X)$, so that a
verification theorem could be proved. Unfortunately we have not been
able to obtain such result. We shall prove though that $V$ is locally
Lipschitz continuous in the $X$-valued variable.

\begin{coroll}
  Under the hypotheses of Proposition
  \ref{prop:Vconv}, the value function $V(s,x)$ is locally Lipschitz continuous with
  respect to $x\in X$. Moreover, the subgradient $\partial V(s,x)$
  with respect to $x$ exists for all $x\in X$ and is locally bounded.
\end{coroll}
\begin{proof}
  The first assertion comes from the fact that a concave locally
  bounded function is continuous in the interior of its effective
  domain (see e.g. Theorem 2.1.3 in \cite{barbu-v}) and $V(s,x)$ is
  finite for all $x\in X$.  Corollary 2.4 in \cite{ET} and the fact
  that $D(\phi)^\circ \subset D(\partial \phi)$ for any concave
  function $\phi$, where $A^\circ$ denotes the interior of a set $A$,
  imply that $V(s,x)$ is locally Lipschitz in $x$, so the assertion on
  $\partial V$ follows.
\end{proof}
Since $V(s,x)\equiv V(s,x_0,x_1)$ is a concave function of $x_0$ for
fixed $s$ and $x_1$, one can also say that $V$ is twice differentiable
almost everywhere with respect to $x_0$, as it follows by the
Busemann-Feller theorem. A similar statement is not true regarding
differentiability with respect to $x_1$, as the Alexandrov theorem is in
general no longer true in infinite dimensions.
We now prove that the value function is continuous with respect to
the time variable. It is possible to prove local Lipschitz
continuity of $V(s,\cdot)$ without appealing to concavity, but
assuming local Lipschitz continuity of $\varphi$.
\begin{prop} The value function $V(s,x)$ is continuous in $s$. Moreover, if
$|\varphi_0(x)-\varphi_0(y)| \leq K(1+R)^m|x-y|$ for all
  $|x|$, $|y| \leq R$, then the function $V(s,x)$ is locally Lipschitz
  continuous with respect to $x$. Furthermore, there exists a constant
  $N=N(K,m)$ such that
  \begin{equation}
    \label{eq:Vlip}
  |V(s,x)-V(s,y)| \leq N(1+R)^m |x-y|
  \end{equation}
  for all $|x|$, $|y| \leq R$.
\end{prop}
\begin{proof}
Recalling that the difference of two suprema is less or equal to the
supremum of the difference, we have
\begin{eqnarray*}
|V(s,x) - V(s,y)| &\leq&
\sup_{z\in\mathcal{U}}|J(s,x;z)-J(s,y;z)|\\
&\leq& \sup_{z\in\mathcal{U}} \E|\varphi(Y^{z,s,x}) - \varphi(Y^{z,s,y})|,
\end{eqnarray*}
and, by Cauchy-Schwarz' inequality,
\begin{eqnarray*}
&& \E|\varphi(Y^{z,s,x}(T)) - \varphi(Y^{z,s,y}(T))| \leq \\
&& \qquad K\Big(\E(1 + |Y^{z,s,x}(T)|^m + |Y^{z,s,y}(T)|^m)^2\Big)^{1/2}
   \Big(\E|Y^{z,s,x}(T) - Y^{z,s,y}(T)|^2\Big)^{1/2}.
\end{eqnarray*}
Arguing as in the proof of Proposition \ref{prop:Vfin}, there exists a
constant $N_1=N_1(K,m)$ such that $\E|Y^{z,s,x}(T)|^m \leq N_1(1+|x|^m)$.
Furthermore, a simple calculation reveals that
$|Y^{z,s,x}(t)-Y^{z,s,y}(t)|=|e^{(t-s)A}(x-y)|$, hence
$$
\E|Y^{z,s,x}(T) - Y^{z,s,y}(T)|^2 \leq M_T^2 |x-y|^2,
$$
and the second claim is proved.
Let us now prove that $V(s,x)$ is continuous in $s$ for a fixed $x$.
Let $s_n \uparrow s$ be a given sequence (the case $s_n \downarrow s$
is completely similar).
Bellman's principle yields
\begin{equation}
  \label{eq:bp}
V(s_n,x) = \sup_{z\in\mathcal{U}} \E\left[
\int_{s_n}^s h(z(r))\,dr + V(s,Y^{z,s_n,x}(s)) \right],
\end{equation}
and choosing a $1/n$-optimal strategy $z_n$ in (\ref{eq:bp}),
we have
\begin{equation}
\label{eq:tano}
\begin{array}{rcl}
\ds \limsup_{n\to\infty}
V(s_n,x)-V(s,x) &\leq&
\ds \limsup_{n\to\infty} \E\int_{s_n}^s h(z_n(r))\,dr\\[10pt]
&& \ds + \limsup_{n\to\infty} \E|V(s,Y^{z_n,s_n,x}(s))-V(s,x)|.
\end{array}
\end{equation}
The first term on the right-hand side in (\ref{eq:tano}) is zero
because $h$ is bounded on $U$. Let us show the also the second term on
the right-hand side of (\ref{eq:tano}) is zero: in fact we have
\begin{eqnarray*}
\E|Y^{z_n,s_n,x}(s) - x|^2 &\leq& K \Big(
\E\int_{s_n}^s |e^{(s-r)A}h(z_n(r))|^2\,dr
+ \E\Big|\int_0^{s-s_n} e^{(s-s_n-r)A}G\,dW(r)\Big|^2 \Big) \\
&\leq& K_1(s-s_n) + q_n \to 0
\end{eqnarray*}
as $n\to\infty$, hecause $h$ is bounded on $U$ and the stochastic
convolution is a Gaussian random variable with covariance operator
going to $0$ as $n\to\infty$. Therefore we also have $Y^{z_n,s_n,x}(s)
\to x$ in probability. By (\ref{eq:Vlip}), $V(s,x)$ is
continuous in $x$ uniformly with respect to $s$, hence
$$
V(s,Y^{z_n,s_n,x}(s)) \to V(s,x)
$$
in probability by the continuous mapping theorem. Moreover, recalling
that $|V(s,x)|\leq N(1+|x|^m)$ and $\E\sup_{t\in[0,T]}
|Y(t)|^m<\infty$, we have
\begin{equation}
\label{eq:porcon}
  \E|V(s,Y^{z_n,s_n,x}(s)) - V(s,x)|^m < \infty,
\end{equation}
hence Vitali's theorem implies
$$
\E|V(s,Y^{z_n,s_n,x}(s)) - V(s,x)| \to 0,
$$
hence $V(s_n,x)-V(s,x) \to 0$ as $n\to \infty$. Furthermore, taking
any $z_0 \in U$, we have from (\ref{eq:bp}),
$$
V(s,x)-V(s_n,x) \leq V(s,x) - (s_n-s)h(z_0) - \E V(s,Y^{z_0,s_n,x}),
$$
which goes to zero as $n \to \infty$ by (\ref{eq:porcon}), and the
claim is proved.
\end{proof}
\begin{rmk}
  Notice that, by the Rademacher theorem in infinite dimensions, the
  previous proposition implies that the value function $V(s,x)$ is
  differentiable in a dense subset of $X$. The local Lipschitz
  continuity of $V$ also implies that (Clarke's) generalized gradient
  of $V(s,x)$ with respect to $x$ is defined everywhere on $X$.
\end{rmk}
For the following proposition, which establishes a monotonicity
property of the value function, we need to define the natural ordering
in $X$: we shall write $x^1\geq x^2$ if $x^1_0 \geq x^2_0$ and $x^1_1
\geq x^2_1$ almost everywhere. Similarly, $x^1>x^2$ if the previous
inequalities hold with the strict inequality sign.
\begin{prop}\label{prop:Vincr}
If $a_1 \geq 0$ and $\varphi_0$ is increasing, then the value
  function $V(s,x)$ is increasing with respect to $x$ in the sense just defined.
\end{prop}
\begin{proof}
  The proof is completely analogous to that of proposition
  \ref{prop:V-incr} below if we prove that $A$ generates a positivity
  preserving semigroup. This is indeed the case: in fact, a direct
  calculation shows that $A$ is the adjoint of the operator $\A$
  defined in section \ref{sec:ex}. Since $a_1\geq 0$, $\A$
  generates a positivity preserving semigroup $S(t)$. It is well known
  that $A$ is the generator of the adjoint semigroup $S(t)^*$. Let
  $x$, $y$ be arbitrary positive elements of $X$. Then
$$
0 \leq \cp{S(t)x}{y} = \cp{x}{S(t)^*y}.
$$
By the arbitrariness of $x$ and $y$, $S(t)^*$ is positivity preserving.
\end{proof}

\subsection{Approximating the value function and the optimal strategy}
Let us now consider the Bellman equation on $X$ associated to the problem of
maximizing (\ref{ex1jbis}), which can be written as
\begin{equation}
\label{eq:HJB0}
\left\{\begin{array}{l}
\ds v_t + {1\over 2}\tr (GG^*v_{xx}) + \cp{A x}{v_x} + H_0(v_x) = 0,
\quad 0 \leq t \leq T\\[8pt]
v(T) = \varphi,
\end{array}\right.
\end{equation}
where $H_0(p) = \sup_{z\in U} (\cp{Bz}{p}+h(z))$.

The main problem with (\ref{eq:HJB0}) is that it is not solvable with
any of the techniques currently available, with the possible exception
of the theory of viscosity solutions. In particular, as of now, one
cannot characterize the value function as the (unique) solution, in a
suitable sense, of equation (\ref{eq:HJB0}). 
As a consequence, we cannot obtain an optimal strategy for the
optimization problem at hand. As a (partial) remedy we develop a
method to approximate the value function and to construct suboptimal
feedback strategies that are asymptotically optimal, in the sense of
Proposition \ref{prop:asop} below.  Let us also briefly recall that,
if we know a priori that a smooth solution to the Bellman equation
exists, then we can apply the verification theorem proved in
\cite{levico}, which in turns allows to obtain precise
characterizations of the optimal strategy (see subsection
\ref{subsec:exex} and section \ref{sec:ex}).

\smallskip

Let us begin proving some approximation results for the value
function $V$. Let $\varepsilon\in]0,1]$ and define $G_\varepsilon:
\erre^2\to X$ as
$$
G_\varepsilon = \left[\begin{array}{cc}
                         \sigma_0 & 0 \\
                         0        & \varepsilon b_1
                      \end{array}
                \right].
$$
Let $W_1$ be a standard real Wiener process independent of $W_0$, set
$W=(W_0,W_1)$, and denote by $\tilde{\mathbb{F}}$ the filtration
generated by $W$. Let $\tilde{\mathcal{U}}$ be the set of
$\tilde{\mathbb{F}}$-adapted processes taking values in $U$.
\par\noindent
Consider the following approximating SDE on $X$:
\begin{equation}
\label{eq:eps}
dY(t) = [A Y(t) + Bz(t)]\,dt + G_\varepsilon\,dW(t),
\quad Y(s)=\bar{x},
\end{equation}
where $z\in\tilde{\mathcal{U}}$ and $0 \leq s \leq t \leq T$.
\par\noindent
For a fixed $z$, let $Y$ and $Y_\varepsilon$ be, respectively,
solutions of (\ref{eq:abstract}) and of (\ref{eq:eps}).  Moreover, let
us define $\varphi_\varepsilon(x)=(\varphi_{0\varepsilon}(x),0)$,
$\varphi_{0\varepsilon}=\tilde{\varphi}_{0\varepsilon}\ast\zeta_\varepsilon$,
$\tilde{\varphi}_{0\varepsilon}(x) =
\varphi_0(x)\chi_{[-1/\varepsilon,1/\varepsilon]}(x)$, and
$h_\varepsilon(x)=(-h_{0\varepsilon}(x),0)$,
$h_{0\varepsilon}(x)=h_0\ast\zeta_\varepsilon(x)$. In particular
$\varphi_{0\varepsilon}\in C^2_c(\erre)$, $h_{0\varepsilon}\in
C^2(\erre)$.  Finally, the approximate objective function and value
function are defined as
$$
J_{\varepsilon_1,\varepsilon_2}(s,x;z) =
\E_{s,\bar{x}}^z\Big[\varphi_{\varepsilon_2}(Y_{\varepsilon_1}(T))
+ \int_s^T h_{\varepsilon_2}(z(t))\,dt \Big],
\qquad
V_{\varepsilon_1,\varepsilon_2}(s,x) =
\sup_{z\in\tilde{\mathcal{U}}} J_{\varepsilon_1,\varepsilon_2}(s,x;z).
$$
In the following we shall set
$\varepsilon=(\varepsilon_1,\varepsilon_2)$ and
$\lim_{\varepsilon\to 0}:= \lim_{\varepsilon_2\to
0}\lim_{\varepsilon_1\to 0}$. Moreover, for $R>0$ we set
$C_R=\left\{x\in X: \;|x_0|\le R \right\}$.
\begin{thm}\label{prop:approx}
  One has $V_\varepsilon(s,x) \to V(s,x)$ as $\varepsilon \to 0$
  uniformly over $s\in[0,T]$, $x\in C_R$, for all $R>0$.
\end{thm}
\begin{proof}
Setting $\eta_{\varepsilon_1}(t)=Y_{\varepsilon_1}(t)-Y(t)$, one has
$$
\eta_{\varepsilon_1}(t) = \varepsilon_1 \int_s^t e^{(t-r)A}B_1\,dW_1(r)
= :\varepsilon_1 \eta(T),
$$
with $B_1:\erre\to X$, $B_1:x \to (0,b_1(\cdot)x)$,
and (suppressing the subscripts on the expectation sign for simplicity)
\begin{eqnarray}
|J_\varepsilon(s,x;z)-J(s,x;z)| &\leq&
|\E[\varphi_{\varepsilon_2}(Y_{\varepsilon_1}(T)) - \varphi(Y(T))]|\nonumber\\
&& + \left|\E\Big[\int_s^T h_{\varepsilon_2}(z(t))\,dt
              - \int_s^T h(z(t))\,dt\Big]\right|\nonumber\\
&\leq& |\E[\varphi_{\varepsilon_2}(Y_{\varepsilon_1}(T))
           - \varphi_{\varepsilon_2}(Y(T))]|
+ |\E[\varphi_{\varepsilon_2}(Y(T)) - \varphi(Y(T))]|\nonumber\\
\label{eq:Je}
&& + \left|\E\Big[\int_s^T h_{\varepsilon_2}(z(t))\,dt
              - \int_s^T h(z(t))\,dt\Big]\right|
\end{eqnarray}
Since
$$
\E|Y_{\varepsilon_1}(T)-Y(T)| =
\varepsilon_1 \E\Big|\int_s^t e^{(t-r)A}B_1\,dW_1(r)\Big| \to 0,
$$
then $Y_{\varepsilon_1}(T) \to Y(T)$ in probability uniformly over
$x\in X$ as $\varepsilon_1\to 0$, and by the continuous mapping
theorem $\varphi_{\varepsilon_2}(Y_{\varepsilon_1}(T)) \to
\varphi_{\varepsilon_2}(Y(T))$ in probability. Moreover
$\varphi_{\varepsilon_2}(Y(T)) \to \varphi(Y(T))$ in probability
uniformly over $x\in C_R$, for all $R>0$, as $\varepsilon_2\to 0$
because $\varphi_{0\varepsilon_2}(x)\to\varphi_0(x)$ $dx$-a.e. in
$\erre$.
Let us now prove that $\varphi_{\varepsilon_2}(Y_{\varepsilon_1}(T))$ is
uniformly integrable with respect to $\varepsilon$. First let us
observe, as it is immediate to show, that there exists $\bar{K}$,
independent of $\varepsilon_2$, such that
$\varphi_{\varepsilon_2}(x)\leq\bar{K}(1+|x|)^m$.
Then we can write
\begin{eqnarray}
\sup_{\varepsilon\in]0,1]^2} \E|\varphi_{\varepsilon_2}(Y_{\varepsilon_1}(T))| &\leq&
\sup_{\varepsilon_1\in]0,1]} \bar{K}\E(1+|Y(T)+\eta_{\varepsilon_1}(T)|)^m\nonumber \\
&\leq& K_1 + K_2 \sup_{\varepsilon_1\in]0,1]}
                 (\E|Y(T)|^m + {\varepsilon_1}^m\E|\eta(T)|^m)\nonumber\\
&=& \label{eq:ui1}
 K_1 + K_2 (\E|Y(T)|^m + \E|\eta(T)|^m) < \infty,
\end{eqnarray}
where we used twice the inequality $|x+y|^m\leq 2^m(|x|^m+|y|^m)$ and
Burkholder-Davis-Gundy's inequality.
Furthermore,
\begin{equation}
\label{eq:ui2}
\sup_{\varepsilon\in]0,1]^2} \E[|\varphi_{\varepsilon_2}(Y_{\varepsilon_1}(T))|\chi_A]
\leq \E[(K_1 + K_2(|Y(T)|^m + |\eta(T)|^m))\chi_A] \to 0
\end{equation}
as $\mathbb{P}(A)\to 0$, because $K_1 + K_2(|Y(T)|^m + |\eta(T)|^m$
has finite expectation.  Then (\ref{eq:ui1}) and (\ref{eq:ui2}) imply
that $\varphi_{\varepsilon_2}(Y_{\varepsilon_1}(T))$ is uniformly
integrable (see e.g. \cite{kall}, lemma 3.10), hence
$$
| \E\varphi_{\varepsilon_2}(Y_{\varepsilon_1}(T)) - \E\varphi(Y(T)) | \to 0
$$
as $\varepsilon\to 0$ (see e.g. \cite{kall}, proposition 3.12).
\par\noindent
Similarly, since
$$
|h_{\varepsilon_2}(z(t,\omega)) - h(z(t,\omega))| \leq
K_1 + K_2 |z(t,\omega)|^m
$$
for all $t\in[0,T]$, $\omega\in\Omega$ and $\E\int_0^T |z(t)|^m <
\infty$ (because $U$ is compact), by the dominated convergence theorem
we have
$$
\E\int_0^T |h_{\varepsilon_2}(z(t)) - h(z(t)|\,dt \to 0
$$
as $\varepsilon_2 \to 0$.

In view of (\ref{eq:Je}) we have thus proved that
$|J_\varepsilon(s,x;z)-J(s,x;z)| \to 0$ uniformly over $s\in[0,T]$,
$x\in C_R$, for all $R>0$ and $z\in\tilde{\mathcal{U}}$, hence also
that $V_\varepsilon(s,x) \to V(s,x)$.
\end{proof}
If the cost function $h_0$ is continuous, one can use a different
regularization, without requiring compactness of $U$.
\begin{prop}
  If $h$ is continuous, then the assertion of theorem \ref{prop:approx}
  holds.
\end{prop}
\begin{proof}
  Let $h_{\varepsilon,\delta}$ be the sup-inf convolution of $h$ (in
  the sense of \cite{LasLio}), that is
  $$
  h_{\varepsilon,\delta}(x) = \sup_{z\in X}\inf_{y\in X}
     \Big(\frac{|z-y|^2}{2\varepsilon} - \frac{|z-x|^2}{2\delta} + h(y) \Big), \qquad 0<\delta<\varepsilon.
  $$
  It is known that $h_{\varepsilon,\delta}$ is differentiable with
  continuous derivative, that $\inf_{x \in X} h(x) \le
  h_{\varepsilon,\delta}(x) \le h(x)$ for all $x \in X$ and that
  $\lim_{\varepsilon,\delta \to 0^+}h_{\varepsilon,\delta}(x) \rightarrow h(x)$
uniformly over $x\in C_R$ (see \cite{LasLio}). Setting
$h_{\varepsilon_2}=h_{\varepsilon_2,\varepsilon_2/2}$, proposition
\ref{prop:Vfin} and the dominated convergence theorem yield
  $$
  \E\int_0^T h_{\varepsilon_2}(z(t))\,dt
  \stackrel{\varepsilon_2\downarrow 0}{\longrightarrow}
  \E\int_0^T h(z(t))\,dt,
  $$
  which implies that the third term on the right-hand side of
  (\ref{eq:Je}) converges to 0 as $\varepsilon \to 0$.
\end{proof}
Theorem \ref{prop:approx} (or its variant), together with the
following result, allow one to approximate the value function $V(s,x)$
in terms of the solutions of a sequence of Bellman equations.
\begin{prop}
  Assume that the hypotheses of theorem \ref{prop:approx} are
  verified. Assume moreover that $h_0$ is strictly convex and
  $\varphi_0$ is concave.  Then the approximate value function
  $V_\varepsilon$ is the unique mild solution (in the sense of
  \cite{FT05}) of the Bellman equation
  $$
  v_t + \frac12 \tr(G_{\varepsilon_1} G^*_{\varepsilon_1} v_{xx})
  + \cp{Ax}{v_x} + H_{0\varepsilon_2}(v_x) = 0,
  \qquad v(T)=\varphi_{\varepsilon_2},
  $$
  where $H_{0\varepsilon_2}(p)=\sup_{z\in U}(\cp{Bz}{p}+h_{\varepsilon_2}(z))$.
\end{prop}
\begin{proof}
Setting
$$
\tilde{B}_{\varepsilon_1}:\erre\to\erre^2, \qquad
\tilde{B}_{\varepsilon_1} = \left[\begin{array}{c} b_0/\sigma_0 \\
\varepsilon_1^{-1}
                              \end{array}
                        \right],
$$
the approximating equation (\ref{eq:eps}) can be rewritten as
\begin{equation}
\label{eq:epss}
dY(t) = [AY(t) + G_{\varepsilon_1} \tilde{B}_{\varepsilon_1} z(t)]\,dt
+ G_{\varepsilon_1}\,dW(t),
\quad Y(s)=\bar{x}.
\end{equation}
The state equation (\ref{eq:epss}), hence also (\ref{eq:eps}), is
covered by the FBSDE approach to semilinear PDEs in Hilbert spaces
(see e.g. \cite{FT-corso}). In order to prove the statement, we
shall verify that hypothesis 7.1 in \cite{FT05} holds true. In
particular, $G_{\varepsilon_1}$ is Hilbert-Schmidt because $b_1\in
L^2([-r,0],\erre_+)$; $\varphi_{\varepsilon_2}$ is Lipschitz because
$\varphi_{0\varepsilon_2} \in C^2_c(\erre)$;
$|\tilde{B}_{\varepsilon_1} z|_X$ is bounded for $z\in U$ because
$b_1\in L^2([-r,0],\erre)$ and $U$ is compact; finally, since $h_0$
is proper and positive, it is immediate to find
$\varepsilon_0\in]0,1]$, $C\geq 0$ such that $h_{\varepsilon_2}(x)
\geq -C$ and $\inf_U h_{\varepsilon_2} \leq C$, for all positive
$\varepsilon_2<\varepsilon_0$.

Since $h$ is strictly convex, for a sufficiently small
$\varepsilon_2$ also $h_{\varepsilon_2}$ is strictly convex.
Therefore we have
$$
g_{\varepsilon_2}(p)= \arg\max_{z\in
U}(\cp{Bz}{p}+h_{\varepsilon_2}(z))=
(h_{\varepsilon_2}')^{-1}(B^*p).
$$
The claim now follows from \cite{FT05}, theorem 7.2 provided we
prove that the closed loop equation
\begin{equation}
\label{eq:CLE} dY(t) = [AY(t) + B g(v_x(t,Y(t))]\,dt +
G_{\varepsilon_1}\,dW(t), \quad Y(s)=\bar{x}.
\end{equation}
admits a solution. In fact this follows as in Theorem 7.2 of
\cite{FT1}.
\end{proof}
\begin{rmk}
  In fact the convexity of $U$ implies that $H_{0\varepsilon_2}\in
  C^1(\erre)$, and hence that $V_\varepsilon \in C^{0,1}([0,T],X)$, as
  in corollary \ref{cor:Vdiff} below.
\end{rmk}

The above approximations do not give a way to construct approximately
optimal strategies for the original problem. In fact, it is well known
that the problem of constructing approximately optimal controls from
the knowledge of an approximate value function is very hard, and in
general unsolved. However, it is possible to construct a (suboptimal)
feedback control for which we have some error control, in the sense
defined below. For a map $f:[0,T]\times X \to U$ such that the
equation
$$
dY(t) = AY(t)\,dt + Bf(t,Y(t))\,dt + G\,dW(t),\qquad Y(s)=x,
$$
admits a mild solution $Y(t)$, let us set $u_f(t)=f(t,Y(t))$ and
$V^f(s,x)= J(s,x;u_f)$. Similarly we define $V_\varepsilon^f(s,x)$.
Let us suppose that we can obtain a feedback law $f:[0,T]\times X
\to U$, which is approximately optimal for the regularized problem,
and let us write $V_\varepsilon \approx V_\varepsilon^f$ to mean
that the two values differ by a small constant. Moreover, recall
that $V \approx V_\varepsilon \approx V_\varepsilon^f$.
\begin{prop}
  Let $f(t,x)$ be Lipschitz in $x$ uniformly over $t$. Then
  $V_\varepsilon^f(s,x) \to V^f(s,x)$ as $\varepsilon\to 0$.
\end{prop}
\begin{proof}
Denote by $Y^f$ and $Y^f_\varepsilon$, respectively, the
solutions of the equations
\begin{eqnarray*}
  dY^f(t) &=& AY^f(t)\,dt + Bf(t,Y^f(t))\,dt + G\,dW(t) \\
  dY^f_\varepsilon(t) &=& AY^f_\varepsilon(t)\,dt
                          + Bf(t,Y^f_\varepsilon(t))\,dt
                          + G_\varepsilon\,dW(t),
\end{eqnarray*}
with $Y^f(s)=Y^f_\varepsilon(s)=x$. Let us assume, without loss of
generality, $s=0$.  Let us show that $Y^f_\varepsilon(t) \to Y^f(t)$
in $L^1(\Omega,\mathbb{P})$, hence in probability, for all
$t\in[0,T]$: by variation of constants we have
\begin{eqnarray*}
|Y^f_\varepsilon(t) - Y^f(t)| &\leq&
\int_0^t |e^{(t-s)A}B(f(s,Y^f_\varepsilon(s))-f(s,Y^f(s)))|\,ds\\
&& + \varepsilon \Big|\int_0^t e^{(t-s)A}B_1\,dW_1(s)\Big|\\
&\leq& \int_0^t m(s)|Y^f_\varepsilon(s) - Y^f(s)|\,ds
       + \varepsilon \Big|\int_0^t e^{(t-s)A}B_1\,dW_1(s)\Big|,
\end{eqnarray*}
where $m(s)=|e^{(t-s)A}|\,|B|\,|f|_{\mathrm{Lip}}$. Taking expectation
on both sides and recalling that the stochastic convolution has finite
mean, Gronwall's lemma yields
\begin{equation}
  \label{eq:gro}
\E|Y^f_\varepsilon(t) - Y^f(t)| \leq \varepsilon N e^{\int_0^T m(s)\,ds}
\to 0
\end{equation}
as $\varepsilon \to 0$, hence $Y^f_\varepsilon(t) \to Y^f(t)$ in
probability for all $t\in[0,T]$.
\par\noindent
By the same arguments used in the proof of theorem \ref{prop:approx}
we obtain that
\begin{equation}
  \label{eq:lim1}
\lim_{\varepsilon_2\to 0}\lim_{\varepsilon_1\to 0}
\E\varphi_{\varepsilon_2}(Y^f_{\varepsilon_1}(T)) = \E\varphi(Y^f(T)).
\end{equation}
Similarly,
\begin{eqnarray*}
|h_{\varepsilon_2}(f(t,Y^f_{\varepsilon_1}(t)))-h(f(t,Y^f(t)))|
&\leq&
|h_{\varepsilon_2}(f(t,Y^f_{\varepsilon_1}(t)))
    - h_{\varepsilon_2}(f(t,Y^f(t)))| \\
&& + |h_{\varepsilon_2}(f(t,Y^f(t)))-h(f(t,Y^f(t)))|
\end{eqnarray*}
and
$$
h_{\varepsilon_2}(f(t,Y^f_{\varepsilon_1}(t)))
\to h_{\varepsilon_2}(f(t,Y^f(t)))
$$
in probability as $\varepsilon_1\to 0$ for all $t\in[0,T]$, because
$Y^f_{\varepsilon_1}(t) \to Y^f(t)$ in probability and
$h_{\varepsilon_2}\circ f(t,\cdot)$ is continuous. Furthermore,
$h_{\varepsilon_2}(f(t,Y^f(t))) \to h(f(t,Y^f(t)))$ $\mathbb{P}$-a.s.
for all $t\in[0,T]$ as $\varepsilon_2\to 0$ because $h_{\varepsilon_2}
\to h$ a.e. on $\erre$.  Since
$|h_{\varepsilon_2}(f(t,Y^f_{\varepsilon_1}(t)))| \leq \sup_{x\in
  U}h_{\varepsilon_2}(x) < \infty$, then
$h_{\varepsilon_2}(f(t,Y^f_{\varepsilon_1}(t)))$ is uniformly
integrable with respect to $\varepsilon_1$, hence
\begin{equation}
  \label{eq:lim2}
\E\int_0^T \Big|
h_{\varepsilon_2}(f(t,Y^f_{\varepsilon_1}(t)))                                 - h_{\varepsilon_2}(f(t,Y^f(t)))
\Big|\,dt \to 0
\end{equation}
as $\varepsilon_1\to 0$.
Finally,
\begin{equation}
  \label{eq:lim3}
\E\int_0^T \Big|h_{\varepsilon_2}(f(t,Y^f(t)))\,dt
- h(f(t,Y^f(t)))\Big|\,dt \to 0
\end{equation}
as $\varepsilon_2\to 0$ by the dominated convergence theorem, taking
into account that
$$
\E \int_0^T |h_{\varepsilon_2}(f(t,Y^f(t)))-h(f(t,Y^f(t)))|\,dt
< \infty,
$$
because $f$ is bounded, as follows by the compactness of $U$.  The
claim now follows by (\ref{eq:lim1}), (\ref{eq:lim2}) and
(\ref{eq:lim3}).
\end{proof}
The previous proposition does not obviously allow one to say that
$u(t)=f(t,Y(t))$ is an approximately optimal feedback map for the
original problem, as $f$ itself in general depends on $\varepsilon_1$,
$\varepsilon_2$.
The next proposition gives quantitative estimates on
$|V^f(s,x)-V^f_\varepsilon(s,x)|$.
\begin{prop}     \label{prop:asop}
  Assume that $f(t,x)$ is Lipschitz in $x$ uniformly over $t$,
  and that $\varphi$, $h$ are Lipschitz continuous. Then there exist
  constants $N=N(|f|_\mathrm{Lip})$ and $\delta=\delta(\varepsilon_2)$
  such that
$$
|V^f(s,x)-V^f_\varepsilon(s,x)| \leq N\varepsilon_1 + \delta(\varepsilon_2)
$$
with $\lim_{\varepsilon_2\to 0} \delta(\varepsilon_2) = 0$.
\end{prop}
\begin{proof}
Let us write
\begin{eqnarray*}
|V^f-V^f_\varepsilon| &\leq&
\E| \varphi_{\varepsilon_2}(Y^f_{\varepsilon_1}(T)) - \varphi(Y^f(T))| \\
&& + \E\int_0^T \Big|h_{\varepsilon_2}(f(t,Y^f(t)))\,dt                           - h(f(t,Y^f(t)))\Big|\,dt \\
&=& I_1+I_2,
\end{eqnarray*}
and
\begin{eqnarray*}
I_1 &\leq& |\E[\varphi_{\varepsilon_2}(Y^f_{\varepsilon_1}(T))
           - \varphi_{\varepsilon_2}(Y^f(T))]|
+ |\E[\varphi_{\varepsilon_2}(Y^f(T)) - \varphi(Y^f(T))]|\\
&=& I_{11} + I_{12}.
\end{eqnarray*}
Then
$$
I_{11} \leq |\varphi_{\varepsilon_2}|_{\mathrm{Lip}}
\E|Y^f_{\varepsilon_1}(T) - Y^f(T)| \leq
|\varphi|_{\mathrm{Lip}} \, \varepsilon_1 e^{\int_0^T m(s)\,ds},
$$
where we used the fact that mollification does not increase the
Lipschitz constant.
We also have
\begin{eqnarray*}
I_{12} &\leq& \E\Big[
|\varphi_{0\varepsilon_2}(Y_0^f(T)) - \varphi_0(Y_0^f(T))|;
|Y_0^f(T)| \leq \varepsilon_2^{-1}
\Big]\\
&& + \E\Big[
|\varphi_{0\varepsilon_2}(Y_0^f(T)) - \varphi_0(Y_0^f(T))|;
|Y_0^f(T)| > \varepsilon_2^{-1}
\Big]\\
&\leq& \delta_1(\varepsilon_2) + \delta_2(\varepsilon_2),
\end{eqnarray*}
where
$$
\delta_1(\varepsilon_2) = \sup_{|x|\leq1/\varepsilon_2}
|\varphi_{0\varepsilon_2}(x) - \varphi_0(x)| < \infty,
$$
as $\varphi_{0\varepsilon_2}$ converges to $\varphi_0$ uniformly
on compact sets, and $\delta_2(\varepsilon_2)$ is defined as follows:
there exist $K_1$, $K_2\geq 0$ such that
$$
\E\Big[
|\varphi_{0\varepsilon_2}(Y_0^f(T)) - \varphi_0(Y_0^f(T))|;
|Y_0^f(T)| > \varepsilon_2^{-1}
\Big]
\leq
\E\Big[ K_1+K_2|Y_0^f(T)|^m; |Y_0^f(T)| > \varepsilon_2^{-1} \Big],
$$
and
\begin{eqnarray*}
Y^f(T) &\leq& e^{TA}x
+ \int_0^T e^{(T-t)A}Bf(t,Y^f(t))\,dt
+  \int_0^T e^{(T-t)A}G\,dW(t) \\
&\leq& e^{TA}x
+ \int_0^T e^{(T-t)A}R\,dt
+ \int_0^T e^{(T-t)A}G\,dW(t) =: \mu_2 + Z_1.
\end{eqnarray*}
Similarly,
$$
Y^f(T) \geq e^{TA}x
+ \int_0^T e^{(T-t)A}r\,dt
+ \int_0^T e^{(T-t)A}G\,dW(t) =: \mu_1 + Z_1,
$$
where $\mu_1$, $\mu_2\in X$, $U \subseteq [r,R]$, and
$Z_1$ is a centered $X$-valued Gaussian random variable.
Denoting by $Z$ the $\erre$-valued components of $Z_1$, we have
$$
Y^f_0(T) - Z \leq (\mu_2)_0,
\qquad
Y^f_0(T) - Z \geq (\mu_1)_0,
$$
hence $|Y^f_0(T) - Z| \leq |(\mu_1)_0| \vee |(\mu_2)_0|=:\mu$, or
equivalently $|Y_0^f(T)| \leq \mu + |Z|$.  In particular, $Z$ is a
centered Gaussian random variable. Then
\begin{eqnarray*}
\E\Big[ K_1+K_2|Y_0^f(T)|^m; |Y_0^f(T)| > \varepsilon_2^{-1} \Big]
&\leq& K_1 \mathbb{P}(|Z|+\mu > \varepsilon_2^{-1}) \\
&& + K_2\E\Big[ (|Z|+\mu)^m; |Z|+\mu > \varepsilon_2^{-1} \Big]\\
&=:& \delta_2(\varepsilon_2).
\end{eqnarray*}
Note that $\delta_2(\varepsilon_2) \to 0$ as $\varepsilon_2\to 0$
since $\E(|Z|+\mu)^m < \infty$.

We have
\begin{eqnarray*}
I_2 &\leq& \Big|
h_{\varepsilon_2}(f(\cdot,Y_{\varepsilon_1}))
- h_{\varepsilon_2}(f(\cdot,Y))
\Big|_{L^1_T}
+ \Big|
h_{\varepsilon_2}(f(\cdot,Y))
- h(f(\cdot,Y))
\Big|_{L^1_T}\\
&=& I_{21}+I_{22},
\end{eqnarray*}
where $L^1_T$ stands for $L^1(\Omega\times[0,T],d\mathbb{P}\times
dt)$.  Recalling again that mollification does not increase the
Lipschitz constant, we also have
$$
I_{21} \leq |h|_{\mathrm{Lip}} |f|_{\mathrm{Lip}}
|Y_{\varepsilon_1} - Y|_{L^1_T}
\leq |h|_{\mathrm{Lip}} |f|_{\mathrm{Lip}}
\varepsilon_1 N \int_0^T e^{\int_0^t m(s)\,ds}\,dt.
$$
Finally, using again the uniform convergence on compact sets of
mollified continuous functions,
\[
I_{22} \leq T \delta_3(\varepsilon_2).
\qedhere
\]
\end{proof}

\subsection{An example with explicit solutions}     \label{subsec:exex}
In this subsection we study in detail the optimal advertising problem
with linear reward and quadratic cost. In particular, we shall assume
$h(z)=-\beta z_0^2$ and $\varphi(x) = \gamma x_0$, with $\beta$,
$\gamma>0$. In \cite{levico} we proved that a solution (in integral
sense) of the HJB equation (\ref{eq:HJB0}) is of the type
$$
v(t,x) = \cp{w(t)}{x} + c(t), \qquad t \in [0,T], \; x \in X,
$$
where $w=(w_0,w_1):[0,T]\to X$ and $c:[0,T]\to\erre$ are given by
\begin{equation}\label{eq:HJBf3}
\left\{
\begin{array}{ll}
\ds w_0'(t) + a_0w_0(t) + \int_{-r}^0 a_1(\xi)w_1(t,\xi)\,d\xi = 0, &
t \in [0,T[ \\[8pt]
w_0(T)=\gamma,\\[8pt]
w_1(t,\xi) = w_0(t-\xi)\chi_{[0,T]}(t-\xi),\\[8pt]
\displaystyle c(t) = \int_t^T \frac{(\cp{B}{w(s)}^+)^2}{4\beta}\,ds,
& t \in [0,T].
\end{array}\right.
\end{equation}
Moreover, the optimal strategy is
\begin{equation}
\label{eq:22}
z^*(t) = \frac{\cp{B}{v_x(t)}^+}{2\beta} =
\frac{\cp{B}{w(t)}^+}{2\beta}, \qquad t \in [0,T]
\end{equation}
(see \cite{levico} for more details).

We extend now the analysis of this specific situation. Let us begin
with a rather explicit characterization of the optimal trajectory,
which could be numerically approximated simply by solving a linear ODE
with delay.
In particular, let $w=(w_0,w_1)$ be the solution of (\ref{eq:HJBf3}).
Then, setting $z^*(t)={1\over 2\beta}\cp{B}{w(t)}^+$, the optimal
trajectory is the $\erre$-valued component $Y_0$ of the (mild)
solution of the abstract SDE
$$ dY(t)=[AY(t) + Bz^*(t)]\,dt + G\,dW(t), $$
which is given by
$$
Y(t) = e^{tA}Y(0) + \int_0^t e^{(t-s)A}Bz^*(s)\,ds
+ \int_0^t e^{(t-s)A}G\,dW(s).
$$
In particular $Y$ is a $X$-valued Gaussian process with mean
and covariance operator
$$
\mu_t = e^{tA}Y(0) + \int_0^t e^{(t-s)A}Bz^*(s)\,ds ,
\quad
Q_t = \int_0^t e^{(t-s)A}GG^*e^{(t-s)A^*}\,ds,
$$
respectively.
It follows that $Y_0$ is a Gaussian process itself with mean
\begin{eqnarray}
\E Y_0(t) &=& \cp{\mu_t}{e_1}_X =
\left\langle e^{tA}Y(0),e_1\right\rangle_X
+ \left\langle \int_0^t e^{(t-s)A}Bz^*(s)\,ds,e_1\right\rangle_X \nonumber\\
\label{eq:opttraj}
&=& \left\langle Y(0),e^{tA^*}e_1\right\rangle_X
+ \int_0^t\left\langle Bz^*(s),e^{(t-s)A^*}e_1\right\rangle_X ds,
\end{eqnarray}
where $e_1=(1,0)\in X$.\\
Since $Y(0)$ is given as in Proposition \ref{prop:equiv} and
$Bz^*(\cdot)$ is also easy to compute ($z^*$ is one dimensional and
$B$ is just multiplication by a fixed vector in $X$), we are left with
the problem of computing $e^{tA^*}e_1$. However, as one can prove by a
direct calculation, the semigroup $e^{tA^*}$ is given by
$$
e^{tA^*}(x_0,x_1(\cdot)) = \Big(\phi(t),\phi(t+\xi)|_{\xi \in [-r,0]}\Big),
$$
where $\phi(\cdot)$ solves the linear ODE with delay
\begin{equation}
\label{eq:sgdelay}
\left\{\begin{array}{l}
\ds {d\phi(t)\over dt} =
a_0 \phi(t) + \int_{-r}^0a_1(\xi)\phi(t+\xi)\,d\xi,
\quad 0\leq t\leq T \\[10pt]
\phi(0)=x_0; \quad \phi(\xi)=x_1(\xi) \;\; \forall\xi\in[-r,0].
\end{array}\right.
\end{equation}
Therefore $e^{tA^*}e_1$ is given by $(\phi(t),\phi(t+\xi)|_{\xi \in
  [-r,0]})$, where $\phi$ solves (\ref{eq:sgdelay}) with initial
condition $x_0=1$, $x_1(\cdot)=0$. Such $\phi$ can be computed
numerically by discretizing (\ref{eq:sgdelay}), and then $\E Y_0(t)$
can be obtained by approximating the integrals in (\ref{eq:opttraj})
with finite sums.

Analogously one can write the variance of the optimal trajectory in
such a way that it can be easily approximated by numerical methods. In
particular, one has
\begin{eqnarray}
\mathrm{Var}\,Y_0(t) &=& \cp{Q_te_1}{e_1}
= \left\langle
\int_0^t e^{(t-s)A}GG^*e^{(t-s)A^*}\,ds\, e_1,e_1
\right\rangle \nonumber\\
&=& \int_0^t
\left\langle e^{(t-s)A}GG^*e^{(t-s)A^*} e_1,e_1 \right\rangle ds \nonumber\\
&=& \int_0^t
\left\langle Ge^{(t-s)A^*}e_1,Ge^{(t-s)A^*}e_1 \right\rangle ds \nonumber\\
&=& \int_0^t \left|Ge^{(t-s)A^*}e_1\right|^2\,ds.\label{eq:optvar}
\end{eqnarray}
Setting $\psi(s)=e^{(t-s)A^*}e_1$, which can be approximated as
indicated before, one finally has
$$
\mathrm{Var}\,Y_0(t) = \sigma^2 \int_0^t \psi(s)^2\,ds.
$$

\medskip

One can also perform simple comparative statics on the
value function. For instance we can compute explicitly its sensitivity
with respect to the (maximal) delay $r$:
\begin{eqnarray*}
\frac{\partial V}{\partial r}(t,x;r) &=&
{\partial \over \partial r}\cp{w_1(t)}{x_1}_{L^2([-r,0])} +
\frac{\partial c}{\partial r}(t;r) \\
&=& w_1(t,-r)x_1(-r)
+ {1 \over 2\beta}
\int_t^T \cp{B}{w(s)}
{\partial\over\partial r}\left(\int_{-r}^0 b_1(\xi)w_1(s,\xi)\,d\xi\right)\,ds \\
&=& w_1(t,-r)x_1(-r) + {b_1(-r) \over 2\beta}
\int_t^T \cp{B}{w(s)} w_1(s,-r)\,ds,
\end{eqnarray*}
where we have used the fact that $\cp{B}{w(t)}^+=\cp{B}{w(t)}$.
\par\noindent
Note that in the above expression everything can be computed
explicitly, as soon as we fix the delay kernel $b_1$.
Let us consider, as an example, the special case of
$b_1(\xi)=b_1\chi_{[-r,0]}(\xi)$, where on the right-hand side, with a
slight abuse of notation, $b_1$ is a positive constant. One has
\begin{eqnarray*}
\frac{\partial V}{\partial r}(t,x;r) &=& w_1(t,-r)x_1(-r) +
{b_1 \over 2\beta} \int_t^T w_1(s,-r)
\Big(b_0w_0(s) + b_1\int_{-r}^0w_1(s,\xi)\,d\xi\Big)\,ds.
\end{eqnarray*}
Furthermore, if we consider the special case of delay in the control
only, that is $a_1(\cdot)=0$, we obtain, after some calculations,
\begin{eqnarray*}
\frac{\partial V}{\partial r}(t,x;r) &=&
\gamma e^{a_0(T-t+r)}x_1(-r)
- \frac{b_1}{4\beta a_0}\gamma^2e^{a_0r}%
\Big(b_0-{b_1\over a_0}(1-e^{a_0r})\Big)%
(1-e^{2a_0(T-t)}),
\end{eqnarray*}
for $t \in [r,T]$.


\medskip

In the special case of $a_1(\cdot)\equiv 0$ an explicit solution of
(\ref{eq:HJBf3}) is easily obtained. This solvability in closed form
then ``propagates'' to other quantities of interest. In fact, note
that (\ref{eq:HJBf3}) reduces to
\begin{equation}
\left\{
\begin{array}{l}
\ds w_0'(t) + a_0w_0(t) = 0 \\[8pt]
w_0(T)=\gamma,
\end{array}\right.
\end{equation}
yielding
$$ w_0(t) = \gamma e^{(T-t)a_0}, $$
and therefore
$$
w_1(t,\xi) = \gamma e^{(T-(t+\xi))a_0}\,\chi_{[0,T]}(t+\xi)
\quad
c(t) = \int_t^T {\cp{B}{w(s)}^2 \over 4\beta}\,ds.
$$
That is, the last three formulae explicitly give a solution of the
HJB equation (\ref{eq:HJB0}) in our specific case.

As a consequence we can also determine the unique optimal feedback control
$z^*$ as follows:
\begin{equation}
\label{eq:27}
z^*(t) = {\cp{B}{v_x}^+ \over 2\beta} = {\cp{B}{w(t)}^+ \over 2\beta}
= {\gamma e^{(T-t)a_0} \over 2\beta} \left[%
b_0 + \int_{-r}^0 b_1(\xi)e^{-a_0\xi}\chi_{[0,T]}(t+\xi)\,d\xi \right].
\end{equation}

The optimal trajectory can be characterized in a completely similar
way as above, with the difference that now we can explicitly write:
$$
e^{tA^*}e_1 = (e^{a_0t},e^{a_0t+\xi}|_{\xi\in[-r,0]}),
$$
hence simplifying (\ref{eq:opttraj}) in the present case. Even
simpler is the expression for the variance of the optimal trajectory,
which can be obtained by (\ref{eq:optvar}):
$$
\mathrm{Var}\,Y_0(t) = \sigma^2 \int_0^t e^{2a_0(t-s)}\,ds
=\frac{\sigma^2}{2a_0} (e^{2a_0t}-1).
$$

\begin{figure}[t]
\includegraphics[width=\textwidth]{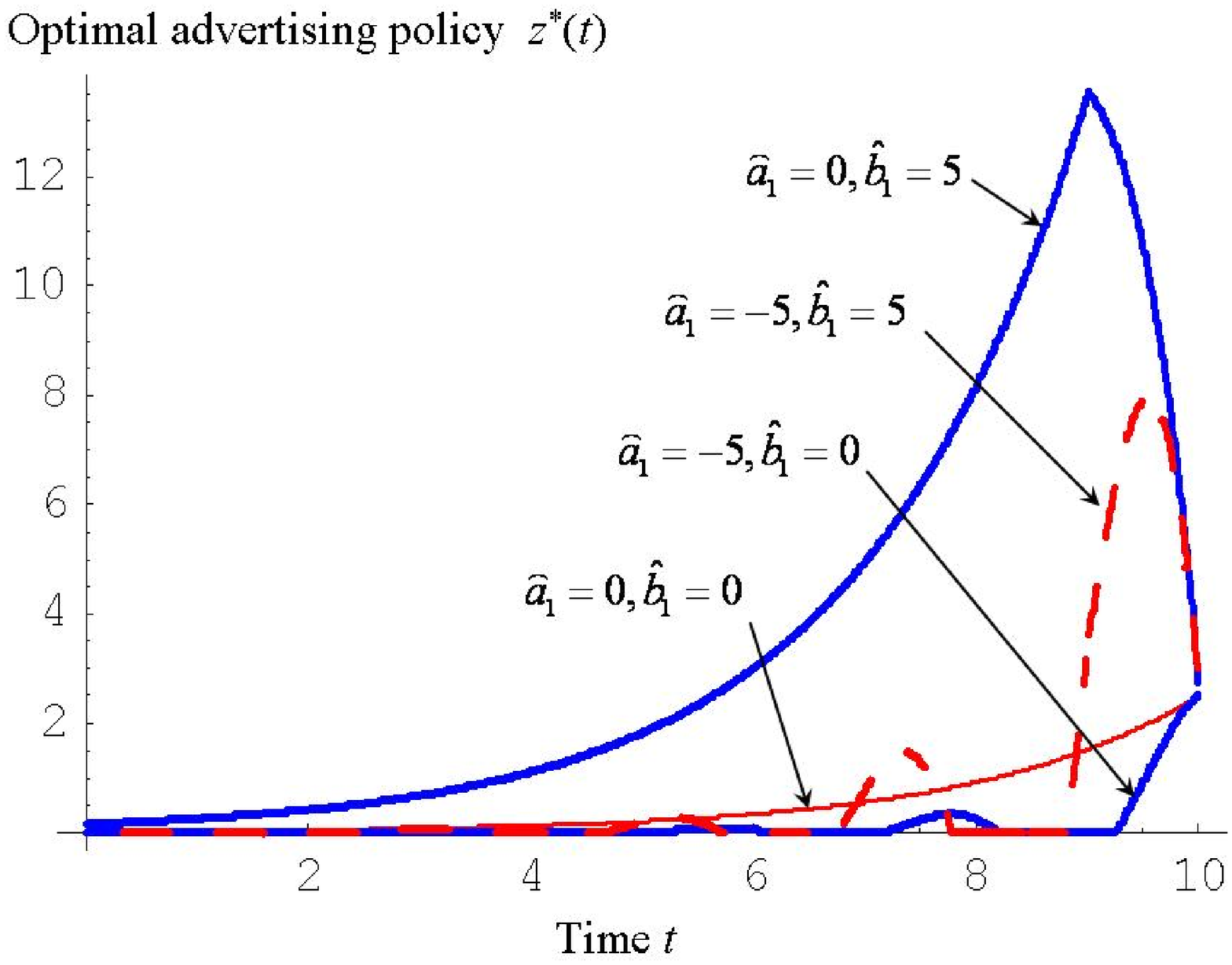}
\caption{Optimal advertising policy in four different ``churn'' settings}
\label{fig:1}
\end{figure}

\begin{figure}[t]
\includegraphics[width=\textwidth]{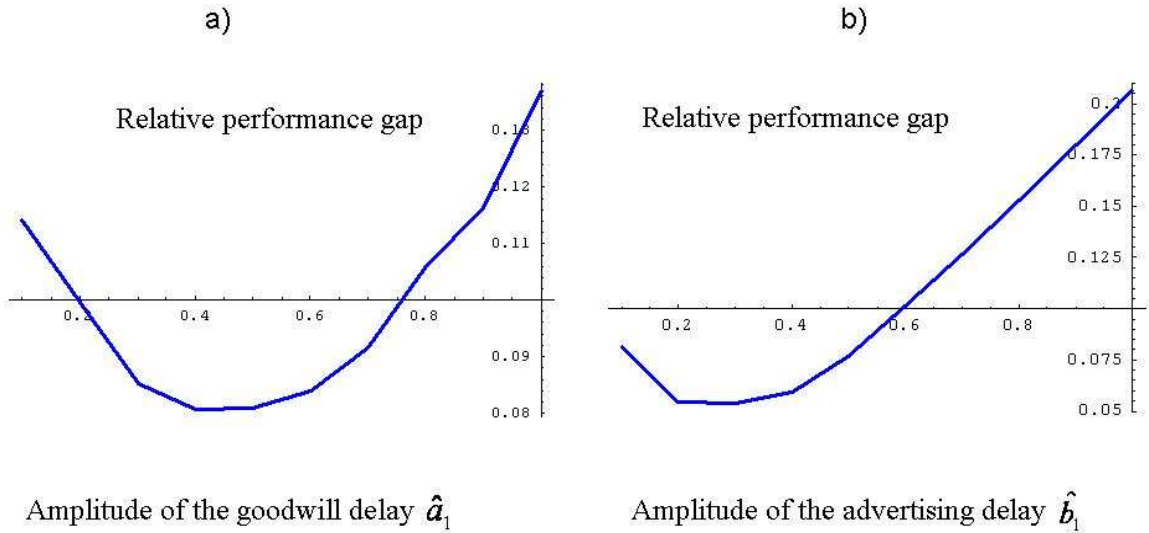}
\caption{Relative performance gap of the ``memoryless'' advertising policy}
\label{fig:2}
\end{figure}

Sharp characterizations of the optimal advertising trajectory as
well as the resulting expected profit functions in the case of
linear reward and quadratic cost function allow for interesting
observations regarding the importance of the proper accounting for
the memory effects in planning the advertising campaign. Figure
\ref{fig:1} displays the optimal advertising spending rates
$z^*(t)$, as expressed by (\ref{eq:22}), with
$a_1(\xi)=\hat{a}_1e^{-|\xi|/\delta_a}$ and
$b_1(\xi)=\hat{b}_1e^{-|\xi|/\delta_b}$ in four different settings:
$a_1 = b_1 = 0$ (``no churn''), $a_1=-5$, $b_1=0$ (``goodwill
churn''), $a_1=0$, $b_1=5$ (``advertising churn''), $a_1=-5$,
$b_1=5$ (``goodwill-advertising churn''). Note that in the absence
of churn, the optimal advertising trajectory, as implied by
(\ref{eq:27}), is a monotone function of time with $z^*(t)=\gamma
b_0/(2\beta)$. While the advertising rates are similar in all
settings in the beginning of the pre-launch period as well as right
before the product launch time $T$, the details of advertising
policies differ dramatically in the middle of the pre-launch period.
For example, in the presence of a strong ``goodwill churn'' the
optimal advertising trajectory takes a characteristic ``impulse''
shape, while in the strong ``advertising churn'' setting the optimal
advertising spending quickly builds up a strong goodwill level in
the middle of the pre-launch region, slowing down significantly
right before the product launch. When the presence of both types of
``churn'' is pronounced, the optimal advertising policy is
represented by a set of advertising sprees with rapidly growing
intensity. Figure \ref{fig:2} illustrates how the strong influence
of memory effects on the shape of optimal advertising policies
translates into performance differences between the optimal
advertising policies and the policies which neglect the presence of
advertising delays. In this figure we plot the relative difference
$$
\frac{V(0,x)-V^0(0,x)}{V(0,x)}
$$
between the optimal expected profit function $V(0,x)$ and the expected
profit value $V^0(0,x)$ obtained by applying, for $t\in [0,T]$, the
advertising policy $z^0(t)=\gamma b_0e^{(T-t)a_0}/(2\beta)$ optimal in
the absence of memory effects (i.e. in the setting $a_1=b_1=0$).  This
relative difference is plotted as a function of the amplitude of the
``goodwill churn'' term $a_1$ (Figure \ref{fig:2}a), and as a function
of the amplitude of the ``advertising churn'' term $b_1$ (Figure
\ref{fig:2}b). The initial goodwill conditions were selected as
$x_0=10$ and $x_1(\xi)= x_0 e^{-|\xi|}$ for $\xi \in [-r,0]$ and the
advertising history $z_0(\xi)$ was set equal to $0$ for
$\xi\in[-r,0]$. We observe that the relative loss of efficiency
associated with the use of the ``memoryless'' policy $z^0(t)$ can be
quite significant -- in the examples we use it exceeds 5\% and can be
as high as 20\% in settings with strong ``churn'' effects.

\section{The case of delay in the state term only}
\label{sec:ex}

In this section we consider a model for the dynamics of goodwill with
forgetting, but without lags in the effect of advertising expenditure
(carryover), i.e. with $b_1(\cdot)=0$ in (\ref{eq:SDDE}). An analysis
of this model was sketched in \cite{levico}, where only an abstract
existence result was given. Here we present a more refined result (see
theorem \ref{thm:L2} below) and obtain some qualitative properties of
the value function, together with a characterization of the optimal
strategies in terms of the value function in two specific cases. In
particular, theorem \ref{thm:L2} formulates sufficient conditions
ensuring that the value function is the unique solution (in a suitable
sense) of the associated Bellman equation and that the optimal
advertising policy is of the feedback type. Let us recall once again
that such a situation is not possible in the more general case
discussed in the previous section. Moreover, in the case of linear
cost function, the optimal control takes a particularly simple
``bang-bang'' form (Corollary \ref{cor:bang}).

Let us also mention that for stochastic control problem with delay
terms in the state variable one can apply both the approach of
Hamilton-Jacobi equations in $L^2$ spaces developed by Goldys and
Gozzi \cite{GolGoz}, and the forward-backward SDE approach of
Fuhrman and Tessitore \cite{FT1}.  We follow here the first approach,
showing that both the value function and the optimal advertising
policy can be characterized in terms of the solution of a Bellman
equation in infinite dimensions.

We assume, for the sake of simplicity, that the goodwill evolves
according to the following equation, where the distribution of the
forgetting factor is concentrated on a point:
\begin{equation}
\label{eq:ex1}
\left\{\begin{array}{l}
dy(t) = \ds [a_0 y(t) + a_1y(t-r)
    + b_0 z(t)]\,dt + \sigma\,dW_0(t), \quad 0\leq s\leq t\leq T \\[10pt]
y(s)=x_0; \quad y(s+\xi)=x_1(\xi), \;\; \xi\in[-r,0].
\end{array}\right.
\end{equation}
The following standard infinite dimensional Markovian reformulation of
this dynamics will turn out to be useful.

Let us define the operator $\A:D(\A)\subset X \to X$ as
$$
\A:(x_0,x_1(\cdot)) \mapsto (a_0x_0+a_1x_1(-r),x_1'(\cdot)),
\qquad
D(\A) = \erre\times W^{1,2}([-r,0];\erre).
$$
It is well-known (see e.g. \cite{DZ96}) that $\A$ is the generator of
a strongly continuous semigroup $S(t)$ on $X$. More precisely, one has
$$
S(t)(x_0,x_1) = (u(t),u(t+\xi)|_{\xi\in[-r,0]}),
$$
where $u(\cdot)$ is the solution of the deterministic delay equation
\begin{equation}
\label{eq:sgA}
\left\{\begin{array}{l}
\ds {du(t)\over dt} = a_0 u(t) + a_1 u(t-r), \quad 0\leq t\leq T \\[10pt]
u(0)=x_0; \quad u(\xi)=x_1(\xi), \;\; \xi\in[-r,0].
\end{array}\right.
\end{equation}
Furthermore, set $\bar{z}=(\sigma^{-1}b_0 z,z_1(\cdot))$, with
$z_1(\cdot)$ a fictitious control taking values in
$L^2([-r,0],\erre)$, and
define $G:X \to X$ as
$$
G: (x_0,x_1(\cdot)) \mapsto (\sigma x_0,0).
$$
Let $W_1$ be a cylindrical Wiener process taking
values in $L^2([-r,0],\erre)$, so that $W=(W_0,W_1)$ is an $X$-valued
cylindrical Wiener process.

Chojnowska-Michalik \cite{choj78} proved the following equivalence result.
\begin{lemma}
\label{lem:equiv}
  Let $Y=(Y_0,Y_1)$ be the unique mild solution of the following stochastic
  evolution equation on $X$:
\begin{equation}
\label{eq:ex1a}
\left\{\begin{array}{l}
dY(t) = (\A Y(t) + G\bar{z}(t))\,dt + G\,dW(t),
              \quad 0\leq s \leq t \leq T,\\[10pt]
Y(s) = x.
\end{array}\right.
\end{equation}
Then $Y_0(t)$ solves the stochastic delay equation (\ref{eq:ex1}).
\end{lemma}
Define $h:X\to\erre$ and $\varphi:X\to\erre$ as
\begin{eqnarray*}
h(x_0,x_1) &=& -h_0(\sigma b_0^{-1}x_0) \\
\varphi(x_0,x_1) &=& \varphi_0(x_0).
\end{eqnarray*}
Then we have, thanks to lemma \ref{lem:equiv},
$$
J(s,x;z) =
\E_{s,x}\left[\varphi(Y(T))+\int_t^T h(\bar{z}(s))\,ds \right],
$$
and
\begin{equation}
\label{ex1j}
V(s,x) = \sup_{\bar{z}\in\mathcal{Z}} J(s,x;\bar{z}),
\end{equation}
where $\mathcal{Z}$ denotes the set of all strategies $z:[0,T]\to
\tilde{U} \times L^2([-r,0],\erre)$ adapted to the filtration
generated by $Y$, and $\tilde{U}$ is the image of $U$ under the action
of the map $x\mapsto \sigma^{-1}b_0x$.

\smallskip

We can now prove some qualitative properties of the value function.
\begin{prop}\label{prop:V-conv}
  If $\varphi_0$ is concave and $h_0$ is convex, then the value
  function $V(s,x)$ is concave with respect to $x$.
\end{prop}
\begin{proof}
Identical to the proof of proposition \ref{prop:Vconv}, thus omitted.
\end{proof}
In the following proposition we use the ordering in $X$ defined right before
Proposition \ref{prop:Vincr}.
\begin{prop}\label{prop:V-incr}
  Let $a_1 \geq 0$ and $\varphi_0$ be increasing. Then the value
  function $V(s,x)$ is increasing with respect to $x$. Moreover, if
  $\varphi_0$ is strictly increasing, then the value function $V(s,x)$
  is strictly increasing with respect to $x$, and $V(s,x^1)=V(s,x^2)$
  if and only if $x^1=x^2$.
\end{prop}
\begin{proof}
Let $x^1\geq x^2$ in the sense just defined. One has
\begin{eqnarray*}
J(s,x^1;z)-J(s,x^2;z) &=&
\E\Big[ \varphi(y^{s,x^1,z}(T)) - \varphi(y^{s,x^2,z}(T)) \Big] \\
&=& \E\Big[ \varphi(e^{\A(T-s)}x^1 + \zeta) -
            \varphi(e^{\A(T-s)}x^2 + \zeta)
      \Big],
\end{eqnarray*}
where $\zeta:=\int_s^T e^{\A(T-t)}G\bar{z}(t)\,dt+%
\int_s^T e^{\A(T-t)}G\,dW(t)$. The assumption $a_1\geq 0$
together with (\ref{eq:sgA}) implies that the semigroup generated by
$\A$ is positivity preserving, i.e.  $x^1\geq x^2$ implies
$e^{\A(T-s)}x^1\geq e^{\A(T-s)}x^2$.  Therefore, by the monotonicity of
$\varphi_0$, one also has $\varphi(e^{\A(T-s)}x^1 + \zeta) \geq
\varphi(e^{\A(T-s)}x^2 + \zeta)$ a.s., hence $J(s,x^1;z) \geq
J(s,x^2;z)$, and finally $V(s,x^1)=\sup_{z\in\mathcal{Z}} J(s,x^1;z)
\geq \sup_{z\in\mathcal{Z}} J(s,x^2;z)=V(s,x^2)$. The other assertions
follow analogously, using (\ref{eq:sgA}).
\end{proof}
\begin{rmk}
  In the above proof the positivity preserving property of the
  semigroup $e^{t\A}$ is crucial, and the assumption $a_1\geq 0$ is
  ``sharp'' in the following sense: if $a_1<0$, one can find $x>0$
  such that $e^{t\A}$ \emph{inverts} the sign, i.e. $e^{t\A}x<0$.

  Moreover, under the assumptions of the theorem, the value function
  is increasing with respect to the real valued component of the
  initial datum. By this we mean that given $x^1 \geq x^2$ with
  $x^1_0>x^2_0$ and $x^1_1=x^2_1$ a.e., then $V(s,x^1)>V(s,x^2)$.
  Therefore one also has $D^-_0V\geq 0$. The subdifferential can be
  replaced by the derivative if we can guarantee that $V$ is
  continuously differentiable with respect to $x_0$. Conditions for
  the continuous differentiability of $V$ with respect to $x$ are
  given in proposition \ref{prop:Vdiff} below.
\end{rmk}

In contrast with the general case considered in the previous section,
if delay enters only the state term, then it is possible to uniquely
solve the associated Bellman equation, and thus to characterize the
value function and construct optimal strategies.  The following
result, which relies on \cite{GolGoz}, gives precise conditions for
the above assertions to hold.
\begin{thm}\label{thm:L2}
  Assume that $h_0$ is convex, let $H_0:\erre\ni p \mapsto
  \sup_{u\in U}(pb_0u-h_0(u))$, and suppose that $a_0 < -a_1 <
  \sqrt{\gamma^2+a_0^2}$, where $\gamma\,\mathrm{coth}\,\gamma = a_0$,
  $\gamma\in]0,\pi[$. Then the value function $V(s,x)$ coincides
  $\mu$-a.e. with the mild solution in $L^2(X,\mu)$ (in the sense of
  \cite{GolGoz}) of the equation
\begin{equation}
\label{eq:ex1-hjb-expl}
\left\{\begin{array}{l}
\ds \partial_t v
+ {1\over 2}\sigma^2 \partial_0^2 v
+ (a_0x_0 + a_1x_1(-r)) \partial_0 v
+ \int_{-r}^0 x'_1(\xi) \partial_1 v(\xi)\,d\xi
+ H_0(\partial_0 v) = 0 \\[10pt]
v(T,x_0,x_1) = \varphi_0(x_0),
\end{array}\right.
\end{equation}
where $\mu$ is a measure of full support on $X$.
Moreover, the optimal strategy admits the feedback representation
\begin{equation}
\label{ex1-zopt}
z^*(t) \in \sigma b_0^{-1}D^-_0H( \sigma\partial_0 v(t,Y^*_0(t),Y_1^*(t) ),
\end{equation}
with $H$ defined in (\ref{eq:h}), provided there exists a solution
$Y_0^*(t)$, $Y_1^*(t)$ of the closed-loop differential inclusion
\begin{equation}
\label{eq:di}
\left\{\begin{array}{l}
\ds dY_0(t) \in \left[a_0Y_0(t) + a_1Y_0(t-r) +
\sigma D^-_0H(\sigma \partial_0 v(t,Y^*_0(t),Y_1^*(t))
\right]dt + \sigma\,dW_0(t) \\[10pt]
\ds dY_1(t)(\xi) = \frac{d}{d\xi}Y_1(t)(\xi).
\end{array}\right.
\end{equation}
\end{thm}
\begin{proof}
By the usual heuristic application of the dynamic programming principle
one can associate to the control problem (\ref{ex1j}) the following
Hamilton-Jacobi-Bellman equation on $X$:
\begin{equation}
\label{ex1hjb}
\left\{\begin{array}{ll}
\ds v_t + {1\over 2}\tr(GG^*v_{xx}) + \cp{\A x}{v_x} + H_0(v_x) = 0,
\quad 0 \leq t \leq T, \\[10pt]
v(T,x) = \varphi(x),
\end{array}\right.
\end{equation}
which coincides, after some calculations, with (\ref{eq:ex1-hjb-expl}).
Note that the Hamiltonian $H_0$ can be regarded as a function of $X$
in $\erre$ and can be equivalently written as
\begin{eqnarray}
H_0(p) = H_0(p_0) &=&
\sup_{z\in \tilde{U}\times L^2([-r,0],\erre)}
\Big(\cp{Gz}{p} + h(z)\Big) \nonumber \\
&=& \sup_{z_0\in \tilde{U}}
\Big(\sigma z_0p_0 - h_0(\sigma b_0^{-1}z_0)\Big).
\label{eq:h0}
\end{eqnarray}
In order to apply the results of \cite{GolGoz}, we also need to define
\begin{equation}
\label{eq:h}
H(q)=H(q_0)=H_0(\sigma^{-1}q_0) =
\sup_{z_0\in\tilde{U}} \Big(z_0q_0 - h_0(\sigma b_0^{-1}z_0)\Big).
\end{equation}
Since $h_0$ is bounded from below, (\ref{eq:h}) implies that $H$ is
Lipschitz continuous (in $\erre$ and in $X$). Moreover, the assumption
on $a_0$, $a_1$ and assumption (i) of section \ref{sec:formul} imply
that the uncontrolled version of (\ref{eq:ex1a}), i.e.
\begin{equation}
\label{eq:ex1auc}
dY(t) = AY(t)\,dt + G\,dW(t),
\end{equation}
admits a unique non-degenerate invariant measure $\mu$ on $X$ (see
\cite{DZ96}), which is Gaussian with mean zero and covariance
operator $Q_\infty=\int_0^\infty e^{sA}GG^*e^{s^*A}\,ds$. In
particular, the restriction of $\mu$ on the $\erre$-valued component
of $X$ has a density
$\rho(x)=\frac{1}{\nu\sqrt{2\pi}}e^{-|x|^2/2\nu^2}$ for some $\nu>0$.
This implies that $\varphi\in L^2(X,\mu)$: in fact,
$$
\int_X |\varphi(x)|^2\,\mu(dx) =
\int_\erre |\varphi_0(x)|^2\rho(x)\,dx \leq
K \int_\erre (1+|x|)^m e^{-|x|^2/2\nu^2}\,dx < \infty.
$$
Therefore, theorems 3.7 and 5.7 of
\cite{GolGoz} yield the existence and uniqueness of a solution in
$L^2(X,\mu)$ of (\ref{ex1hjb}), or equivalently of
(\ref{eq:ex1-hjb-expl}), which coincides $\mu$-a.e. with the value
function $V$.
Finally, observing that the maximum in (\ref{eq:h}) is reached by
$D^-_0H(q_0)$ (setting, if needed, $h_0(x)=+\infty$ for $x\not\in U$),
a slight modification of the proof of theorem 5.7 in \cite{GolGoz}
shows that the optimal strategy is given by $\bar{z}_0^*(t) \in
D^-_0H(\sigma\partial_0 v(t,Y^*_0(t),Y_1^*(t))$, where
$Y^*=(Y_0^*,Y_1^*)$ is a solution (if any) of the stochastic
differential inclusion (\ref{eq:di}). The relation $z^*(t)=\sigma
b_0^{-1}\bar{z}_0^*(t)$ thus completes the proof.
\end{proof}
Let us briefly comment on the previous result: the HJB equation
(\ref{eq:ex1-hjb-expl}) is ``genuinely'' infinite dimensional, i.e. it
reduces to a finite dimensional one only in very special cases. For
example, by the results in \cite{lari}, (\ref{eq:ex1-hjb-expl})
reduces to a finite dimensional PDE if and only if $a_0=-a_1$.
However, under this assumption, we cannot guarantee the existence of a
non-degenerate invariant measure for the Ornstein-Uhlenbeck semigroup
associated to (\ref{eq:ex1auc}). Even more extreme would be the
situation of distributed forgetting time: in this case the HJB
equation is finite dimensional only if the term accounting for
distributed forgetting vanishes altogether.
Moreover, note that if $a_1$ is negative, i.e. it can be interpreted
as a deterioration factor, the assumption of the theorem says that
$a_1$ cannot be ``much more negative'' than $a_0$. On the other hand,
if $a_1$ is positive, then the improvement effect as measured by $a_1$
cannot exceed the deterioration effect as measured by $|a_0|$.  In
essence, the condition on $a_0$, $a_1$, which is needed to ensure
existence of an invariant measure for equation (\ref{eq:ex1auc}), does
not impose severe restrictions on the dynamics of goodwill.

\smallskip

If the data of the problem are smoother, a different approach allows
one to obtain regularity of the value function.
\begin{prop}\label{prop:Vdiff}
  Assume that $\varphi_0\in C^1(\erre)$,
  $|\varphi_0'(x)|\leq K(1+|x|)^m$, and $H_0\in C^1(\erre)$.
  Then $V\in C^{0,1}([0,T]\times X)$.
\end{prop}
\begin{proof}
  Follows by the regularity results for solutions of semilinear
  partial differential equations in Hilbert spaces obtained through
  the FBSDE approach. In particular, denoting by $C$ a positive
  constant, boundedness of $U$ implies that $|b_0 z|<C$, $|h_0(z)|
  \leq K(1+|z|)^m < C$, and finally $H_0$ is Lipschitz as follows by
  \begin{eqnarray*}
    |H_0(p)-H_0(q)| &=& |\sup_{z\in U}(\cp{Bp}{z}-h_0(z))
                       - \sup_{z\in U}(\cp{Bp}{z}-h_0(z))| \\
    &\leq& \sup_{z\in U} |\cp{(p-q)}{B^*z}| \leq C|p-q|.
  \end{eqnarray*}
  Then all hypotheses of \cite{FT-corso}, theorem 4.3.1, are
  satisfied, which yields the claim.
\end{proof}
%
\begin{coroll}\label{cor:Vdiff}
  Let $\varphi_0$ be as in proposition \ref{prop:Vdiff} and $h_0$
  strictly convex. Then $V\in C^{0,1}([0,T]\times X)$.
\end{coroll}
\begin{proof}
  Since convexity implies continuity in the interior of the domain,
  then $h_0(U)$ is bounded. Extending $h_0$ as $h_0(x)=+\infty$ for
  $x\not\in U$, $h_0$ is clearly $1$-coercive, hence $H$ is convex and
  finite on the whole $\erre$ (\cite{lema}, prop. E.1.3.8). The strict
  convexity of $h_0$ implies that $H_0$ is continuously differentiable
  in the interior of its domain, i.e. on $\erre$ (\cite{lema}, thm.
  E.4.1.1). Then the smoothness of $V$ follows again by
  \cite{FT-corso}, theorem 4.3.1.
\end{proof}
\begin{rmk}
  We should also mention that in the framework of the FBSDE approach
  to HJB equations (\cite{FT-corso}), if the Hamiltonian $H$ and the
  terminal condition $\varphi$ satisfy some smoothness and boundedness
  conditions, then we do not need the assumption about the existence
  of an invariant measure for the uncontrolled state equation. The
  approach used above (\cite{GolGoz}), while requiring the
  existence of the above mentioned invariant measure, allows for more
  singular data (for instance one could choose $\varphi_0(x)=-M$,
  $M>\!\!>0$, for $x\in\erre_-$, and $\varphi_0(x)\geq 0$ for
  $x\in\erre_+$).
\end{rmk}

\medskip

In general, obtaining explicit expressions of the value function $V$
trying to solve (\ref{eq:ex1-hjb-expl}) is impossible. However, under
specific assumptions on the model we can obtain stronger
characterizations, at least from a qualitative point of view, of the
value function and/or of the optimal strategy.
\begin{coroll}
Assume that $h_0(x)=\beta x^2$ and $U=[0,R]$, $R<\infty$. Then
the optimal strategy is given by
\begin{equation}
\label{eq:L2-ex1}
z^*(t) =
\left\{\begin{array}{ll}
\ds 0, & D_0V^* < 0 \\[8pt]
\ds \frac{b_0D_0V^*}{2\beta}, & 0 \leq D_0V^* \leq
      2b_0^{-1}\beta R \\[8pt]
\ds R, & D_0V^* > 2b_0^{-1}\beta R,
\end{array}\right.
\end{equation}
where $V^*:=V(t,Y_0^*(t),Y_1^*(t))$.
\end{coroll}
\begin{proof}
One has
\begin{eqnarray*}
H_0(p) &=& \sup_{0\leq z_0\leq \tilde{R}} \Big(\cp{Gp}{z}+h(z)\Big)
= \sup_{0\leq z_0 \leq \tilde{R}} (\sigma p_0 z_0 - \tilde\beta z_0^2) \\
&=& \ds \frac{(\sigma p_0)^2}{4\tilde\beta}
        I_{\{0 \leq p_0 \leq \frac{2\tilde\beta\tilde{R}}{\sigma} \}}
     + (\sigma p_0 \tilde{R} - \tilde\beta\tilde{R}^2)
        I_{\{ p_0 > {2\tilde\beta\tilde{R} \over \sigma} \}}
\end{eqnarray*}
where $\tilde{R}=\sigma^{-1}b_0R$ and $\tilde\beta:=\sigma^2 b_0^{-2}\beta$.
Therefore
$$
H(q) = q_0^2/4\tilde\beta \, I_{\{0 \leq q_0 \leq 2\tilde\beta\tilde{R} \}}
       + (q_0 \tilde{R} - \tilde\beta\tilde{R}^2) \, I_{\{q_0 > 2\tilde\beta\tilde{R}\}}
$$
and
$$
DH(q) = q_0/2\tilde\beta I_{\{0 \leq q_0 \leq 2\tilde\beta\tilde{R}\}}
+ \tilde{R} I_{\{ q_0 > 2\tilde\beta\tilde{R}\}}.
$$
Theorem \ref{thm:L2} now yields (\ref{eq:L2-ex1}).
\end{proof}
Note that whenever $\varphi_0$ is increasing, we get $D^-_0V^*\geq 0$, hence
the optimal control is either linear in $D_0V^*$ or constant for
$D_0V^*$ over a threshold.

\begin{coroll}     \label{cor:bang}
Assume that $h_0(x)=\beta x$. Then the optimal strategy is of
the bang-bang type and is given by
$$
z^*(t) =
\left\{\begin{array}{ll}
\ds 0, & D_0V^* < \sigma b_0^{-2}\beta \\[6pt]
\ds \rho, & D_0V^* = \sigma b_0^{-2}\beta \\[6pt]
\ds R, & D_0V^* > \sigma b_0^{-2}\beta,
\end{array}\right.
$$
where $\rho$ is an arbitrary real number.
\end{coroll}
\begin{proof}
Setting $\tilde{R}=\sigma^{-1}b_0R$ and
$\tilde\beta:=\sigma^2 b_0^{-2}\beta$, one has
$$
H_0(p) = \sup_{0\leq z_0 \leq \tilde{R}} (\sigma p_0 z_0 - \tilde\beta z_0)
= (\sigma p_0-\tilde\beta)\tilde{R} \, I_{\{p_0 > \tilde\beta/\sigma\}}
$$
and $H(q) = (q_0-\tilde\beta)\tilde{R} \, I_{\{q_0 > \tilde\beta\}}$,
thus
$$
D^-H(q) = \erre \, I_{\{q_0 = \tilde\beta\}}
          + \tilde{R}\, I_{\{q_0 > \tilde\beta\}}.
$$
Theorem \ref{thm:L2} yields the conclusion.
\end{proof}

In general, even specifying a functional form of $h_0$, an explicit
solution of the HJB equation for arbitrary $\varphi_0$ is not
available, hence the above expressions of the optimal control strategy
in terms of the value function and their corresponding qualitative
properties are the ``best'' that one can expect, at least in the cases
we have considered.

\section{Concluding remarks}

A number of deterministic advertising models allowing for delay effects have
been proposed in the literature. However, the corresponding problems in the
stochastic setting have not been investigated. One of the reasons is
certainly that a theory of continuous-time stochastic control with delays
has only been developed recently, following two approaches. The first
approach is based on the solution of an associated infinite-dimensional
Hamilton-Jacobi-Bellman equation in spaces of integrable functions (see
\cite{GolGoz}). The other one relies on the analysis of an appropriate
infinite-dimensional forward-backward stochastic differential equation (see
\cite{FT1}). Both approaches, however, cannot be applied to problems with
distributed lag in the effect of advertising.

Problems with memory effects in both the state and the control have
been studied first by Vinter and Kwong \cite{VK} (in a deterministic
LQ setting), and by Gozzi and Marinelli \cite{levico} (in the case of
linear stochastic dynamics and general objective function). A general
theory of solvability of corresponding HJB equations is currently not
available, while an infinite-dimensional Markovian reformulation and a
``smooth'' verification theorem have been proved in \cite{levico}.

In this paper we have concentrated on deriving qualitative properties
of the value function (such as convexity, monotonicity with respect to
initial conditions, smoothness). For specific choices of the reward
and cost functions, we obtain more explicit characterizations of value
function and optimal state-control pair.

While our work makes a substantial initial step in the analysis of the
stochastic advertising problems with delays, more remains to be done.
Potential extensions of the present work include the analysis of problems
with budget constraints as well as problems of advertising through multiple
media outlets with different delay characteristics.

\bibliographystyle{amsplain}
\bibliography{ref}

\def\polhk#1{\setbox0=\hbox{#1}{\ooalign{\hidewidth
  \lower1.5ex\hbox{`}\hidewidth\crcr\unhbox0}}}
\providecommand{\bysame}{\leavevmode\hbox to3em{\hrulefill}\thinspace}
\providecommand{\MR}{\relax\ifhmode\unskip\space\fi MR }
\providecommand{\MRhref}[2]{%
  \href{http://www.ams.org/mathscinet-getitem?mr=#1}{#2}
}
\providecommand{\href}[2]{#2}
\begin{thebibliography}{10}

\bibitem{overad1}
D.~Aaker and J.~M. Carman, \emph{Are you overadvertising?}, Journal of
  Advertising Res. \textbf{22} (1982), 57--70.

\bibitem{barbu-v}
V.~Barbu and Th. Precupanu, \emph{Convexity and optimization in {B}anach
  spaces}, second ed., D. Reidel, Dordrecht, 1986. \MR{87k:49045}

\bibitem{bass2}
F.~M. Bass, \emph{Optimal advertising expenditure implications of
  simultaneous-equation regression analysis}, Oper. Res. \textbf{19} (1971),
  822--831.

\bibitem{bassclark}
F.~M. Bass and D.~G. Clark, \emph{Testing distributed lag models of advertising
  effect}, Journal of Market. Res. \textbf{9} (1972), 298--308.

\bibitem{bassparsons}
F.~M. Bass and L.~J. Parsons, \emph{Simultaneous-equation regression analysis
  of sales and advertising}, Applied Economics \textbf{1} (1969), 103--124.

\bibitem{choj78}
A.~Chojnowska-Michalik, \emph{Representation theorem for general stochastic
  delay equations}, Bull. Acad. Polon. Sci. S\'er. Sci. Math. Astronom. Phys.
  \textbf{26} (1978), no.~7, 635--642.

\bibitem{DP-K}
G.~Da~Prato, \emph{Kolmogorov equations for stochastic {PDE}s}, Birkh\"auser
  Verlag, Basel, 2004. \MR{MR2111320 (2005m:60002)}

\bibitem{DZ96}
G.~Da~Prato and J.~Zabczyk, \emph{Ergodicity for infinite-dimensional systems},
  Cambridge University Press, Cambridge, 1996. \MR{MR1417491 (97k:60165)}

\bibitem{ET}
I.~Ekeland and R.~T{\'e}mam, \emph{Convex analysis and variational problems},
  SIAM, Philadelphia, PA, 1999. \MR{2000j:49001}

\bibitem{elsanosi}
I.~Elsanosi, \emph{Stochastic control for systems with memory}, Dr. Scient.
  thesis, University of Oslo, 2000.

\bibitem{sethi}
G.~Feichtinger, R.~Hartl, and S.~Sethi, \emph{Dynamical {O}ptimal {C}ontrol
  {M}odels in {A}dvertising: {R}ecent {D}evelopments}, Management Sci.
  \textbf{40} (1994), 195--226.

\bibitem{FT1}
M.~Fuhrman and G.~Tessitore, \emph{Nonlinear {K}olmogorov equations in infinite
  dimensional spaces: the backward stochastic differential equations approach
  and applications to optimal control}, Ann. Probab. \textbf{30} (2002), no.~3,
  1397--1465. \MR{MR1920272 (2003d:60131)}

\bibitem{FT-corso}
\bysame, \emph{Backward stochastic differential equations in finite and
  infinite dimensions}, Lecture Notes, Politecnico di Milano, 2004.

\bibitem{FT05}
\bysame, \emph{Generalized directional gradients, backward stochastic
  differential equations and mild solutions of semilinear parabolic equations},
  Appl. Math. Optim. \textbf{51} (2005), 279--332.

\bibitem{GolGoz}
B.~Goldys and F.~Gozzi, \emph{Second order parabolic
  {H}amilton-{J}acobi-{B}ellman equations in {H}ilbert spaces and stochastic
  control: {$L^2_\mu$} approach}, Stoch. Processes Appl. \textbf{116} (2006),
  no.~12, 1932--1963.

\bibitem{levico}
F.~Gozzi and C.~Marinelli, \emph{Stochastic optimal control of delay equations
  arising in advertising models}, Stochastic partial differential equations and
  applications (G.~Da~Prato and L.~Tubaro, eds.), Marcel Dekker, 2005.

\bibitem{griliches}
Z.~Griliches, \emph{Distributed lags: A survey}, Econometrica \textbf{35}
  (1969), 16--49.

\bibitem{hartl}
R.~F. Hartl, \emph{Optimal dynamic advertising policies for hereditary
  processes}, J. Optim. Theory Appl. \textbf{43} (1984), no.~1, 51--72.
  \MR{85h:90042}

\bibitem{HS}
R.~F. Hartl and S.~P. Sethi, \emph{Optimal control of a class of systems with
  continuous lags: dynamic programming approach and economic interpretations},
  J. Optim. Theory Appl. \textbf{43} (1984), no.~1, 73--88. \MR{85g:49024}

\bibitem{overad4}
J.~D. Herrington and W.~A. Dempsey, \emph{Comparing the current effects and
  carryover of national-, regional-, and local-sponsor advertising}, Journal of
  Advertising Research \textbf{45} (2005), 60---72.

\bibitem{lema}
J.-B. Hiriart-Urruty and C.~Lemar{\'e}chal, \emph{Fundamentals of convex
  analysis}, Springer-Verlag, Berlin, 2001.

\bibitem{kall}
O.~Kallenberg, \emph{Foundations of modern probability}, Probability and its
  Applications (New York), Springer-Verlag, New York, 1997. \MR{MR1464694
  (99e:60001)}

\bibitem{KSlibro}
V.~B. Kolmanovski\u{\i} and L.~E. Sha\u{\i}khet, \emph{Control of systems with
  aftereffect}, AMS, Providence, RI, 1996. \MR{1 415 834}

\bibitem{Koyck}
L.~M. Koyck, \emph{Distributed lags and investment analysis}, North-Holland,
  Amsterdam, 1954.

\bibitem{larssen}
B.~Larssen, \emph{Dynamic programming in stochastic control of systems with
  delay}, Stoch. Stoch. Rep. \textbf{74} (2002), no.~3-4, 651--673.
  \MR{2003h:93080}

\bibitem{lari}
B.~Larssen and N.~H. Risebro, \emph{When are {HJB}-equations in stochastic
  control of delay systems finite dimensional?}, Stochastic Anal. Appl.
  \textbf{21} (2003), no.~3, 643--671. \MR{1 978 238}

\bibitem{LasLio}
J.-M. Lasry and P.-L. Lions, \emph{A remark on regularization in {H}ilbert
  spaces}, Israel J. Math. \textbf{55} (1986), no.~3, 257--266. \MR{MR876394
  (88b:41020)}

\bibitem{leone}
R.~P. Leone, \emph{Generalizing what is known about temporal aggregation and
  advertising carryover}, Market. Sci. \textbf{14} (1995), G141--G150.

\bibitem{overad2}
L.~M. Lodish, M.~Abraham, S.~Kalmenson, J.~Livelsberger, B.~Lubetkin,
  B.~Richardson, and M.~E. Stevens, \emph{How tv advertising works: A
  meta-analysis of 389 real world split cable tv advertising experiments},
  Journal of Marketing Res. \textbf{32} (1995), 125--139.

\bibitem{overad3}
X.~Luo and N.~Donthu, \emph{Benchmarking advertising efficiency}, Journal of
  Advertising Research \textbf{41} (2001), 7--18.

\bibitem{admarket}
{TNS} {M}edia {I}ntelligence, \emph{{TNS} {M}edia {I}ntelligence reports {U.S.}
  advertising expenditures increased by 3.0 percent in 2005},
  http://www.tns-mi.com/news/02282006.htm (March 1, 2006).

\bibitem{NA}
M.~Nerlove and J.~K. Arrow, \emph{Optimal advertising policy under dynamic
  conditions}, Economica \textbf{29} (1962), 129--142.

\bibitem{sethi-VWcomp}
A.~Prasad and S.~P. Sethi, \emph{Competitive advertising under uncertainty: a
  stochastic differential game approach}, J. Optim. Theory Appl. \textbf{123}
  (2004), no.~1, 163--185. \MR{MR2100268}

\bibitem{VW}
M.~L. Vidale and H.~B. Wolfe, \emph{An operations-research study of sales
  response to advertising}, Operations Res. \textbf{5} (1957), 370--381.
  \MR{19,514d}

\bibitem{VK}
R.~B. Vinter and R.~H. Kwong, \emph{The infinite time quadratic control problem
  for linear systems with state and control delays: an evolution equation
  approach}, SIAM J. Control Optim. \textbf{19} (1981), no.~1, 139--153.

\end{thebibliography}

\end{document}